\documentclass{amsart}
\input xy
\xyoption{all}
\newcommand{\Ps}{\mathbf{P}}
\newcommand{\Q}{\mathbf{Q}}
\newcommand{\C}{\mathbf{C}}
\newcommand{\HH}{\mathbf{H}}
\newcommand{\Z}{\mathbf{Z}}
\newcommand{\F}{\mathbf{F}}
\newcommand{\kv}{\mathbf{k}}
\newcommand{\mv}{\mathbf{m}}
\newcommand{\av}{\mathbf{a}}

\newcommand{\bv}{\mathbf{b}}

\renewcommand{\phi}{\varphi}
    \newtheorem{Lem}{Lemma}[section]
    \newtheorem{Prop}[Lem]{Proposition}
    \newtheorem{Thm}[Lem]{Theorem}
    
    \newtheorem{Cor}[Lem]{Corollary}

   \theoremstyle{definition}
    \newtheorem{Def}[Lem]{Definition}
    \newtheorem{Not}[Lem]{Notation}
    \newtheorem{Exa}[Lem]{Example}
    \newtheorem{Rem}[Lem]{Remark}

\DeclareMathOperator{\rig}{rig}

\DeclareMathOperator{\Gal}{Gal}
\DeclareMathOperator{\spec}{Spec}
\DeclareMathOperator{\spa}{span}
\DeclareMathOperator{\red}{red}
\DeclareMathOperator{\sing}{sing}
\DeclareMathOperator{\Frob}{Frob}

\DeclareMathOperator{\trace}{trace}
\DeclareMathOperator{\dR}{dR}

\begin{document}
\title[Zeta-function of deformations of Fermat hypersurfaces]{The zeta-function of monomial deformations of Fermat hypersurfaces}
\date{\today}
\author{Remke Kloosterman}
\address{Institut f\"ur Algebraische Geometrie, Universit\"at Hannover, Welfengarten 1, D-30167, Hannover, Germany}
\email{kloosterman@math.uni-hannover.de}

\begin{abstract}
This paper intends to give a mathematical explanation for  results on the zeta-function of some  families of varieties recently obtained in the context of Mirror Symmetry \cite{Cand2}, \cite{Kad}. In doing so, we obtain concrete and explicit  examples for some results recently used in algorithms to count points on smooth hypersurfaces in $\Ps^n$. 

In particular, we extend the monomial-motive correspondence of Kadir and Yui and we give explicit solutions to the $p$-adic Picard-Fuchs equation associated with monomial deformations of Fermat hypersurfaces.

As a by-product (Theorem~\ref{thmPoincare}) we obtain Poincar\'e Duality for the rigid cohomology of certain singular affine varieties.
\end{abstract}
\thanks{This work was partially supported by the DFG Schwerpunktprogramm ``Globale Methoden
in der komplexen Geometrie'' under grant HU 337/5-3.
The author would like to thank Klaus Hulek, Shabnam Kadir, Orsola Tommasi, Jaap Top and an anonymous referee for their remarks and suggestions for improvements on a previous version of this paper.}
\maketitle
\section{Introduction}\label{SecInt}

One of the families under consideration in this paper is the famous one-parameter family (Dwork family) of quintic threefolds $X_{\overline{\lambda}} \subset \Ps^4_{\F_q}$ given by 
\begin{equation}\label{maineq}  x_0^5+x_1^5+x_2^5+x_3^5+x_4^5+\overline{\lambda} x_0x_1x_2x_3x_4 =0,\end{equation}
where $\overline{\lambda}\in \F_q$ is a parameter.
Candelas et al. \cite{Cand2} observed that the zeta-function of this variety can be written as
\[ \frac{R_1(t,\overline{\lambda})R_2(t,\overline{\lambda})^{20}R_3(t,\overline{\lambda})^{30}}{(1-t)(1-qt)(1-q^2t)(1-q^3t)},\]
where the $R_i$ are of degree 4. Candelas et al. gave expressions in $\overline{\lambda}$ for the zeroes of the $R_i$:
to explain this, note that we can lift this family to a family over the ring $\Z_q$ of Witt vector over $\F_q$. 
 This enables us to consider this family as a family in $\Ps^4$ over the field of fractions $\Q_q$ of $\Z_q$. 
 Assume that  $\overline{\lambda}\in \F_q$ is chosen such that $X_{\overline{\lambda}}$ is smooth. 
Denote by $\lambda$  the Teichm\"uller lift of $\overline{\lambda}$. Specifically, Candelas et al. show that the zeroes of the zeta-function of $X_{\overline{\lambda}}$ can be expressed in certain solutions of the $p$-adic Picard-Fuchs equation (associated with the family $X_{\overline{\lambda}}$) evaluated at  $\lambda$.
 
 This  fact was proved in a more general context, but less explicitly, by N. Katz \cite{Katz}. 
His description of the zeta-function in terms of the Picard-Fuchs equation is exploited by Lauder \cite{LauDef}
  in order to give an algorithm to count points on smooth hypersurfaces in $\Ps^n$.

Some other families are investigated by Kadir \cite{Kad}. She obtained similar results. From this, one might conjecture that various factors of the zeta-function are enumerated  by so-called (admissible) monomial types modulo certain equivalence relations. We come back to this in Subsection~{\ref{ssFac}}.

Kadir and Yui \cite{YuiKad} noticed that monomial types are occurring in the study of several objects related to (\ref{maineq}), for example in the Picard-Fuchs equation or in the enumeration of the factors of the zeta-function. In the case  $\lambda=0$, they also appear in the enumeration of the Jacobi sums  needed to compute  the number    of points of the variety at $\lambda=0$. They proved a certain correspondence between these monomial types for Fermat varieties. Our aim is to present a different view on the above mentioned phenomena. 

We should mention that N.~Katz \cite{KatzAno} and Rojas-Leon and Wan \cite{RojWan} studied the zeta-function of  families similar to (\ref{maineq}) by using ($\ell$-adic) hypergeometric sheaves. We recommend \cite{KatzAno} for a discussion on previous results on the Dwork family.

The main object of study in this paper are families $X_{\overline{\lambda}}/\F_q$ defined by the vanishing of polynomials of the form
\begin{equation}\label{mainform} F_{\overline{\lambda}} :=\sum_{i=0}^n x_i^{d_i} + \overline{\lambda} \prod_i x_i^{a_i} \end{equation}
in a weighted projective space $\Ps:=\Ps(w_0,\dots,w_n)$, with $w_id_i=d$ for all $i$, the $a_i$ are non-negative and $\sum w_ia_i=d$; moreover, we assume that $\gcd(q,d)=1$. Such families will be called one-parameter monomial deformations of a Fermat hypersurface.
For the rest of the introduction fix such a weighted projective space, and such a one-parameter deformation of a Fermat hypersurface. Let $\mathbf{a}$ denote the vector $(w_0a_0,w_1a_1,\dots,w_na_n)\in (\Z/d\Z)^{n+1}$. We call $\mathbf{a}$ the \emph{deformation vector}.

The main technical result of this paper implies that the $p$-adic Picard-Fuchs equation associated with such a family is a generalized hypergeometric differential equation. We refer to Subsection~\ref{ssDef} and Section~\ref{SecComp} for  more on this.

Let $U_{\overline{\lambda}}:=\Ps\setminus X_{\overline{\lambda}}$. Since
\[ Z(X_{\overline{\lambda}},t)Z(U_{\overline{\lambda}},t)=Z(\Ps,t),\]
we have that $Z(X_{\overline{\lambda}},t)$ is uniquely determined by $Z(U_{\overline{\lambda}},t)$. Hence from now on we will only discuss how to calculate $Z(U_{\overline{\lambda}},t)$.

\subsection{Choice of the cohomology theory} \label{ssCoh}
The Lefschetz fixed point formula allows us to  prove statements on the zeta-function by considering the action of geometric Frobenius on certain cohomology groups. 
Very often one uses  \'etale cohomology. This is particularly useful when one wants to compare results in characteristic $p>0$ with results in characteristic 0, or if one wants to consider Galois-representations on  certain $\ell$-adic vector spaces.

However, for our purposes it seems more natural to use $p$-adic cohomology theories instead. One can represent  cohomology classes of a variety over a finite field $\F_q$ by differential forms with coefficients in $\Q_q$. This allows us to perform several (basic) analytic tricks when computing with cohomology classes.

To be more precise, let $\lambda$ be a lift of $\overline{\lambda}$ to $\Q_q$, let $F_\lambda$ be a lift of $F_{\overline{\lambda}}$ and $U_{\overline{\lambda}}$.  
Since  $U_\lambda$ is affine,  we can define Monsky-Washnitzer cohomology (cf. Section~\ref{SecMW}) groups $H^i(U_{\lambda},\Q_q)$. The elements in $H^i(U_{\lambda},\Q_q)$ are differential forms with $\Q_q$-coefficients. There is a lift $\Frob_q$ of Frobenius acting on these groups.

To illustrate how explicitly one can compute with Monsky-Washnitzer cohomology, we proceed to produce a basis for $H^i(U_\lambda)$.
Let $\Omega:=\left( \prod_j{x_j}\right)\sum (-1)^i w_i \frac{dx_0}{x_0}\wedge \frac{dx_1}{x_1}\wedge \dots \wedge \widehat{\frac{dx_i}{x_i}} \wedge \dots \wedge \frac{dx_n}{x_n} $. 
\begin{Prop}\label{BasisProp} Suppose $X_{\overline{\lambda}}$ is quasi-smooth. Then the cohomology groups $H^i(U_\lambda,\Q_q)$  are zero except for $i=0,n$. The group $H^0(U_\lambda,\Q_q)$ is one-dimensional and Frobenius acts trivially on it. The following set is a basis for $H^n(U_\lambda,\Q_q)$:
\[ \left\{ \frac{\prod_{i=0}^n x_i^{k_i}}{(F_\lambda)^t}\Omega \colon 0\leq k_i< d_i-1 \; \forall i,\; \sum_i w_i(k_i+1) =td \right\} .\]
\end{Prop}
This basis will be called the {\em standard basis}. We are not aware of a proper reference for this standard fact in our context. We prove this Proposition in Section~\ref{SecMW}.
Proposition~\ref{BasisProp} is a combination of Theorem~\ref{thmLefschetzHyperplane} and  Proposition~\ref{prpBasis}.

The proof  is based on the fact that for quasi-smooth $X_{\overline{\lambda}}$ we have that de Rham cohomology of $U_\lambda$ with $\Q_p$ coefficients is isomorphic to the Monsky-Washnitzer cohomology of $U_{\overline{\lambda}}$ \cite{BC}. By a theorem of Steenbrink \cite{Ste} we have the following isomorphism
 \[ H^n_{\dR}(U_\lambda)\cong  \oplus_{t>0} H^0(\Omega^n(tX_\lambda ))/d H^0(\Omega^{n-1}((t-1)X_\lambda) ).\]
 The vector space on the right-hand-side is very well understood. 
 
However, if $X_{\overline{\lambda}}$ is not quasi-smooth then the dimension of the right-hand-side depends on the choice of the lift $\lambda$. If we choose $\lambda$ in such a way that $X_\lambda$ is not quasi-smooth then the right-hand-side is infinite-dimensional.
In that case one needs to add more relations to get an isomorphism with $H^n(U_\lambda)$. Which relations one needs to add is not very well understood.

A vector \[\kv:=(\overline{w_0(k_0+1)},\dots,\overline{w_n(k_n+1)}) \in \prod_i (w_i\Z/d\Z)\]
is called an  admissible monomial type if for all $i$, we have $k_i\not \equiv -1 \bmod d_i$ and $\sum \overline{w_i(k_i+1)} \equiv 0 \bmod d$. Fix an admissible monomial type $\kv$. Take elements $k_i \in \Z$ satisfying $0\leq k_i \leq  d_i-2$ and $ k_i\equiv \overline{k_i} \bmod d_i$. Then
with $\kv$ we associate the standard basis vector
\[ \omega_\kv :=\frac{\prod x_i^{k_i}}{(F_\lambda)^t}\Omega.\]

\begin{Rem}The results mentioned in Section~\ref{SecMW} imply that
\[ Z(U_{\overline{\lambda}},t)=\frac{\left(\det\left(I-q^n(\Frob_q^*)^{-1}t \mid H^n(U_\lambda,\Q_q)\right) \right)^{(-1)^{n+1}}}{(1-q^n t)} .\]
From here on   we formulate our results in terms of the characteristic polynomial of $q^n(\Frob_q^*)^{-1}$ on $H^n(U_\lambda,\Q_q)$, rather than in terms of $Z(U_{\overline{\lambda}},t)$.
\end{Rem}

\subsection{Deformation behavior}\label{ssDef}
We produce a solution to the $p$-adic Picard-Fuchs equation that turns out to give us a description of the dependence of $\lambda$ of the action of Frobenius on $H^n(U_\lambda)$, where $\lambda$ is in the $p$-adic unit disc.

Following N. Katz \cite{Katz}, we consider the commutative diagram
\[\xymatrix {H^n(U_{\lambda^q}) \ar[r]^{\Frob_q^*} \ar[d]_{A(\lambda^q)}& H^n(U_\lambda) \ar[d]_{A(\lambda)} \\ H^n(U_0) \ar[r]^{\Frob_q^*} & H^n(U_0)}\]
where $\lambda$ is on a small $p$-adic disc around the origin, and $A$ is a solution to the Picard-Fuchs equation associated with the family $X_\lambda$. Using $p$-adic analytic continuation we can extend $A(\lambda)^{-1}\Frob_{q,0}^*A(\lambda^q)$ to the closed unit disc, although $A(\lambda)$ itself cannot be extended to the $p$-adic unit disc.

Let $\lambda_0\in \Q_q$ be the Teichm\"uller lift of some element $\overline{\lambda_0}\in\F_q$. Then $\lambda_0^q=\lambda_0$, hence the above diagram implies that the action of $\Frob_q$ on $H^n(U_{\lambda_0})$ can  be recovered from the $p$-adic analytic continuation of $A(\lambda)^{-1}\Frob_{q,0}^*A(\lambda^q)$. Therefore, to determine the zeta-function of $X_{\overline{\lambda_0}}$ we need to know the Frobenius action in the Fermat-case (see \ref{ssFac}) and  compute the correct solution of the Picard-Fuchs equation.

We describe the action of $A(\lambda)$ on the standard basis. We call two monomial types $\mathbf{k}$ and $\mathbf{m}$ {\em strongly equivalent} if and only if there is a $j_0$ such that $\mathbf{k}-\mathbf{m}=j_0\mathbf{a}$, where $\av$ is the deformation vector (see above).

\begin{Thm} \label{ThmHG} Let $\kv$ be an admissible monomial type. Write  $A(\lambda)\omega_\kv =\sum c_\mv(\lambda) \omega_\mv$, where  the sum is taken over all admissible monomial types. Then $c_\mv(\lambda)$  is non-zero only if $\mathbf{k}$  and $\mathbf{m}$ are  strongly equivalent. If this is the case then $c_\mv(\lambda)$ is of the form $c_0\lambda^{j_0} F(\alpha_i;\beta_j;\lambda^d c_1^d)$, with $F$ a $p$-adic generalized hypergeometric function with parameters $\alpha_i,\beta_j$ and $j_0\in \{0,1,\dots d-1\}$ is chosen such that $\mathbf k-\mathbf m=\overline{j_0}  \mathbf{a}$.
\end{Thm}

Explicit formulas for the $\alpha_i, \beta_j, c_0$ and $c_1$ are given in Lemma~\ref{redformLem} and Proposition~\ref{PrpCoeff}. See Section~\ref{SecComp} for a proof of Theorem~\ref{ThmHG}.

In our proof  we exploit the fact that there is a  straightforward  way of computing in groups like $H^n(U_\lambda)$, relying on the fact that this group is a quotient of a module of differentials over a {\em power series} ring. This allows us to perform  some easy analytic operations that would be impossible in a module of differentials over a polynomial ring. 

\subsection{Factorization of the zeta-function}\label{ssFac} 
We call the case $\lambda=0$ the Fermat case.
One can show that $\Frob_q^*$ on $H^n(U_0,\Q_q)$ sends  the standard basis vector $\omega_\mathbf{k}$ to a constant $c_{\mathbf{k},q}$ times the standard basis vector $\omega_{\overline{q}\kv}$. Hence, if $q\equiv 1 \bmod d$ then the standard basis is a basis of  eigenvectors for $\Frob_q^*$.  In this case Theorem~\ref{ThmHG}  tells us that for every admissible monomial type $\kv$ the operator $\Frob_{q,\lambda}$ fixes the subspace spanned by the $\omega_\mv$, where $\mv$ is strongly equivalent to $\kv$.

The general case is slightly different, for this we introduce another equivalence relation: we call two monomial types $\mathbf{k}$ and $\mathbf{m}$ {\em weakly equivalent} if $j_0 \in \Z/d\Z$ and invertible $s,t\in (\Z/d\Z)^*$ exists such that $s\mathbf{k}+t\mathbf{m}=j_0\mathbf{a}$.   

\begin{Thm}\label{thmWE} Let $\kv$ be an admissible monomial type. Write  \[\Frob_{q,\lambda}\omega_\kv =\sum c_\mv(\lambda) \omega_\mv,\] where  the sum is taken over all admissible monomial types. Then $c_\mv(\lambda)$  is non-zero only if $\mathbf{k}$  and $\mathbf{m}$ are  weakly equivalent. \end{Thm}
This is a weak form of Theorem~\ref{PrpFactor}.
Theorem~\ref{thmWE} implies that the zeta-function of $U_\lambda$ can be factored  (as a rational function with $\Q(\lambda)$-coefficients) in such a way that each factor corresponds to a weak-equivalence class. If one only considers the zeta-function over fields containing all $d$-th roots of unity, then there is a factorization of the zeta-function of $U_\lambda$ such that each factor corresponds to a strong-equivalence class.

Explicitly determining the constants $c_{\mathbf{k},q}$ is actually very hard. In some cases it is known that the eigenvalues of $\Frob_q^*$ correspond to the Fourier coefficients of a modular form. For example if $n=2$, $w=(1,1,1)$ and $d=3$, then $X_0$ is the $j=0$ elliptic curve $x_0^3+x_1^3+x_2^3$. Also the case $n=3$, $w=(1,1,1,1)$, $d=4$ and the case $n=5$, $w=(1,1,1,1)$, $d=3$ are known to correspond to modular forms (see \cite{KH}, \cite{SI}).

A more general result on $c_{\mathbf{k},q}$ is due to Weil:
Assume that $\F_q\supset \F_p(\zeta_d)$. Let $\chi$ be the $d$-th power residue symbol. Let $\kv$ be an admissible monomial type. Let $\kv_i$ be the $i$-th entry of $\kv$, i.e., $w_i(k_i+1)$. Then
\[ J_{\kv,q}:=(-1)^{n+1} \sum_{(v_1,\dots,v_n)\in \F_q^n \colon \sum_i v_i=-1} \chi(v_1)^{\kv_1} \chi(v_2)^{\kv_2}\dots \chi(v_n)^{\kv_n}.\]

 The following Theorem coincides with Corollary~\ref{CorJac}.

\begin{Thm}\label{ThmJac} Assume $q$ is chosen such that $\F_q=\F_p(\zeta_d)$. Let $\mathbf k$ be an admissible monomial type. Let $S$ be the set of monomial types that are weakly equivalent to $\mathbf k$. Then the sets $\{ q^{n-1}/c_{\mathbf m,q} \colon m \in S\}$ and $\{J_{\mathbf m, q}  \colon m \in S\}$ coincide.
\end{Thm}

\subsection{Monomial-motive correspondence}\label{ssMotMon}

We call $\mathbf{b}\in (\Z/d\Z)^{n+1}$ an admissible automorphism type if  $\mathbf{b}=(w_0b_0,w_1b_1,\dots,w_nb_n) \in(\Z/d\Z)^{n+1}$ is such that $\sum w_ib_ia_i \equiv 0 \bmod d$.  Define $\sigma_{\mathbf{b}}$ to be the automorphism
\[ [x_0:x_1:\dots:x_n]\mapsto [\zeta_{d}^{w_0b_0} x_0:\zeta_{d}^{w_1b_1}x_1: \dots:\zeta_d^{w_nb_n}x_n].\]
We call two monomial types $\mathbf{k}$ and $\mathbf{m}$ {\em distinguishable by automorphisms} if there exists an admissible automorphism type $\mathbf{b}\in (\Z/d\Z)^{n+1}$ such that
\[  \sigma_{\mathbf{b}}\left(\prod x_i^{k_i}\right) =\prod x_i^{k_i} \mbox{ and } \sigma_{\mathbf{b}}\left(\prod x_i^{m_i} \right) \neq \prod x_i^{m_i}. \]

\begin{Thm}\label{ThmWeDi} Two monomial types $\mathbf{k}$ and $\mathbf{m}$ are weakly equivalent if and only if $\mathbf{k}$ and $\mathbf{m}$ are not distinguishable by automorphisms.\end{Thm}

This result enables us to give a different proof for the monomial-motive correspondence of Kadir and Yui \cite{YuiKad}, and to generalize it as follows: fix an admissible monomial type $\kv$. Let $G_\kv$ be the group of automorphisms of the form $\sigma_{\mathbf{b}}$ that fix $\omega_\kv$. Then the subspace of $H^n(U)$ fixed by $G_\kv$ is the spanned by  the $\omega_\mv$ such that  $\mv$ is weakly equivalent to $\kv$. 
This can be also extended to the level of motives, i.e., we find a submotive $\mathfrak{h}(U_\lambda/G_\kv)$ of the (Chow-)motive $\mathfrak{h}(U_\lambda)$. Moreover, we obtain that
\[ \mathfrak{h}(U_\lambda) = \oplus_{[\kv]} \mathfrak{h}(U_\lambda/G_\kv), \]
where we sum over all the weak-equivalence classes. 
 
 Kadir and Yui decompose $\mathfrak{h}(U_\lambda/G_{\kv})$ further. To explain this, we need to change our context, and consider our family $X_\lambda$ over the field $\Q$ of rational numbers. Then the Galois group $\Gal(\Q(\zeta_d)/\Q)$ acts non-trivially on $G_\kv$, and this enables us to find correspondences in $CH^n(U_\lambda/G_\kv\times U_\lambda/G_\kv)$ that decompose $\mathfrak{h}(U_\lambda/G_\kv)$ into smaller motives. It is easy to see that each such motive corresponds to a strong-equivalence class of monomial types. This correspondence between admissible monomial types and submotives of $\mathfrak{h}^n(U_\lambda)$ is called by Kadir and Yui  \emph{monomial-motive correspondence}. 
They also relate  monomial types with the Picard-Fuchs equation. For this issue we refer to Subsection~\ref{ssDef}.

Kadir and Yui \cite{YuiKad} could only prove their monomial-motive correspondence if $X_\lambda$ is a Calabi-Yau hypersurface of dimension 3 and $\lambda=0$. The above discussion  extends this correspondence to any quasi-smooth member of a  one-parameter monomial deformation of a Fermat hypersurface in a weighted projective space,  for any degree $d$ such that $w_i|d$ for all $i$ and provided that  the characteristic does not divide $d$. 

Kadir and Yui  prove the monomial-motive correspondence  using Jacobi sums. We take a more direct approach using subgroups of the automorphism group.

This paper is organized as follows: in Section~\ref{SecNot} we fix some notation and list some standard definitions.
In Section~\ref{SecMW} we discuss Monsky-Washnitzer cohomology groups and recall some of the properties of these groups.
In Section~\ref{SecDeform} we recall Katz' result on the deformation of the zeta-function  of a hypersurface in $\Ps^n$.
In Section~\ref{SecComp} we make Katz' result explicit. 
In Section~\ref{SecDia} we discuss the Frobenius action on the cohomology of a Fermat hypersurface and prove some results on the structure of the zeta-function of a monomial deformation of a Fermat hypersurface.

\section{Notation}\label{SecNot}
Fix   once and for all :
\begin{itemize}
\item a prime $p$ (the characteristic) and a positive integer $r$,
\item an integer $n$ (the dimension of the ambient space),
\item a vector $(w_0,w_1,\dots,w_n)\in \Z^{n+1}$ such that none of the $w_i$ is divisible by $p$.
\item an integer $d$ divisible by all the $w_i$ and $p$ does not divide $d$.
\end{itemize}
Set $q=p^r$ and $d_i:=d/w_i$. Let $\Q_q$ denote the unique unramified extension of degree $r$ of $\Q_p$.  Let $w$ denote the total weight, i.e., $w:=\sum w_i$. 
 Let $\Ps_{\F_q}:=\Ps^n_{\F_q}(w_0,\dots,w_n)$ be the associated weighted projective space over the finite field $\F_q$.

\begin{Def}
A {\em monomial type} $\mathbf{m}=(\overline{m_0},\dots,\overline{m_n})$ is an element of $\prod_i w_i\Z/d\Z$ such that   $\sum \overline{m_i}=0$ in $\Z/d\Z$. Choose representatives $m_i\in \Z$ of $\overline{m_i}$ such that $0\leq m_i <d$. The {\em relative degree} of $\mathbf{m}$ is $\sum m_i/d$.
\end{Def}

Fix  once and for all a monomial type $\mathbf{a}$ of relative degree $1$, with at least 2 non-zero entries. We call $\mathbf{a}$ the \emph{deformation vector}. Let $a_i$ be integers such that $0\leq a_i<d_i$ and $\mathbf{a}\equiv (\overline{w_0}\overline{a_0},\dots,\overline{w_n}\overline{a_n})$. Set
\[ F_{\lambda}:= \sum x_i^{d_i}+ \lambda \prod x_j^{a_j}.\]
Let $F:=F_0$. If $\overline{\lambda}\in \F_q$, denote by $X_{\overline{\lambda}}$ the zero set of $F_{\overline{\lambda}}$ in $\Ps$.
If $\lambda\in \Q_q$, denote by $X_\lambda$ the zero set of $F_\lambda$.
  Let $U_{\overline{\lambda}}$ be the complement $\Ps \setminus X_{\overline{\lambda}}$. 
 Let ${U_\lambda}$ be the complement $\Ps \setminus {X_\lambda}$.

Let $\Omega:= \prod_i x_i\sum_j (-1)^j w_j \frac{dx_0}{x_0}\wedge \dots \wedge \widehat{\frac{dx_j}{x_j}} \wedge \frac{dx_n}{x_n}$.

\begin{Def} A monomial type $\mathbf{k}$ is called {\em admissible} if there exist integers $k_i, i=0,\dots n$  such that $0\leq k_i\leq d_i-2$ and $\mathbf{k}:=(w_0(k_0+1),\dots,w_n(k_n+1))$. Let $t$ be the relative degree of $\kv$.
With $\mathbf{k}$ we associate the differential form
\[ \omega_k:=  \frac{\prod x_i^{k_i}}{F_\lambda^t} \Omega.\]\end{Def}

Denote by $(a)_m$ the Pochhammer symbol $a(a+1)\dots(a+m-1)$.

 \begin{Def} \label{defquo}
Let $\pi: \Ps^n \to \Ps$ be the natural quotient map sending $x_i$ to $x_i^{w_i}$. Let $G:=\times {\mathbf{\mu}}_{w_i}/\Delta$ be the Galois group associated with this quotient. We call $\pi$ the {\em standard quotient map} and $G$ the {\em group associated with the standard quotient map}.
\end{Def} 

\section{Monsky-Washnitzer cohomology}\label{SecMW}
We will not define rigid cohomology in complete detail, but give a simplified presentation for the case of quasi-smooth hypersurfaces. For a good introduction to the theory of rigid cohomology we refer to \cite{BerSMF} and \cite{BerFin}.

 Since $U_\lambda$ is affine, we can write $U_\lambda=\spec R_\lambda$, with 
\[ R_\lambda=\Q_q [\lambda,Y_0,\dots, Y_m]/(G_{1,\lambda},\dots,G_{k,\lambda}).\]

\begin{Def}
Fix $\lambda_0$ in the closed $p$-adic unit disc
set
\[ R^{\dagger}_{\lambda_0} = \frac{\{ H \in \Q_q[[Y_0,\dots,Y_m]] \colon \mbox{ the radius of convergence of } H \mbox{ is at least } r>1 \}}{(G_{1,\lambda_0},\dots,G_{k,\lambda_0})}.\]
 Then $R^{\dagger}_{\lambda_0}$ is called the {\em overconvergent completion} (or weak completion) of $R_{\lambda_0}$.

 Let $\pi$ be the standard quotient map and $G$ its associated group (cf. Definition~\ref{defquo}).
Set $S:=S^{\dagger}_{\lambda_0}$ to be the overconvergent completion of the coordinate ring of $\Ps^n\setminus \pi^{-1}(X_{\lambda_0})$. Then on the module of differential forms $\Omega^i_{S}$ there is a natural $G$-action. Set $\Omega^i_{R}=(\Omega^i_{S})^G$. The {\em  $i$-th Monsky-Washnitzer cohomology group $H^i(U_{\lambda_0},\Q_q)$} is the $i$-th cohomology group of the complex $\Omega^\bullet_{R}$.
\end{Def}

\begin{Not}Let $X \subset \Ps$ be a quasi-projective variety. 
 Denote by $H^i_{\rig}(X)$ the {\em $i$-th rigid cohomology group}  of $X$ and by $H^i_{\rig,c}(X)$ the {\em $i$-th rigid cohomology group with compact support}  of $X$, as defined in \cite{BerSMF}.
\end{Not}

There exists a second, equivalent, definition of $H^i(U_{\lambda_0},\Q_q)$. This goes  as follows: since $U_{\lambda_0,\sing}$ is affine, there is a ring $S$ such that $U_{\lambda_0,\sing}=\spec S$. Let $S^\dagger$ be an overconvergent completion of $S$. Let $\iota:\spec R^\dagger_{\lambda_0}\setminus\spec S^\dagger \to \spec R_{\lambda_0}^\dagger$ be the inclusion. Let $\Omega^i_{\spec R^\dagger_{\lambda_0}}$ be the sheaf $\iota_*\Omega_{\spec R^\dagger_{\lambda_0} \setminus \spec S^\dagger}$. Then define the Monsky-Washnitzer cohomology groups $H^i(U_{\lambda_0},\Q_q)$ as the cohomology groups of the complex obtained by taking global sections. The proof that these two definitions are equivalent is very similar to \cite[2.2.4]{Dol}.

\begin{Def} Let $R$ be a ring over $\Z_q$. Let $\pi$ be the maximal ideal of $\Z_q$. A {\em lift of Frobenius}  is a ring homomorphism 
 $\Frob_q^*: R \rightarrow R$ such that  its reduction modulo $\pi$
 \[ \Frob_q^* \bmod \pi :  R\otimes_{\Z_q} \F_q  \rightarrow R\otimes_{\Z_q} \F_q \]
 is well-defined and equals $x\mapsto x^q$. 
 \end{Def}

Fix a lift of Frobenius $\Frob_q^*$ to $R_{\lambda_0}^\dagger$, such that $\Frob_q^*(\lambda)=\lambda^q$. By abuse of notation we denote by $\Frob^*_q$ also  the induced morphism on $H^i(U_{\lambda_0},\Q_q)$.

 \begin{Prop}\label{prpCompare} There is a natural isomorphism
 \[ H^i_{\rig}(U_{\overline{\lambda_0}},\Q_q)\cong H^i(U_{\lambda_0},\Q_q)\]
 which is compatible with the action of Frobenius.
 \end{Prop}

 \begin{proof}Similar to the proof of \cite[Proposition  1.10]{BerFin}.
 \end{proof}

\begin{Def}
Let $K$ be a field. Let $G\in K[x_0,\dots,x_n]$ be a weighted homogeneous polynomial (with weights $(w_0,\dots,w_n)$). Let $Y$ be the hypersurface $G=0$ in $\Ps$. Then $Y$ is said to be {\em quasi-smooth} if the affine cone $\spec K[X_0,\dots,X_n]/G$ is smooth or has exactly one singular point, namely $(0,0,\dots,0)$.
\end{Def}
\begin{Rem}\label{RemSm}
If $\Ps=\Ps^n$ then a hypersurface $X\subset \Ps$ is quasi-smooth if and only if it is smooth.\end{Rem}

An easy calculation shows
\begin{Lem} \label{lemQS}  Let $I=\{i\in \{0,1,\dots,n\}\colon a_i\not \equiv 0 \bmod p\}$. Let $g=\gcd_{i\in I}(a_iw_i)$ and $d':=d/g$.
If there is a non-zero $a_i$ such that $a_i\equiv 0 \bmod p$ then $X_\lambda$ is quasi-smooth for all $\lambda$. Otherwise,
$X_{\overline{\lambda}}$ is quasi-smooth if and only if 
 \[ \overline{\lambda}^{d'}\neq  \frac{(-1)^{d'} \overline{d}^{d'}}{\prod_{i\in I} (\overline{a_iw_i})^{a_iw_i/g}} .\] \end{Lem}

 \begin{proof} Consider the partial derivative of $F$ with respect to $x_j$. If $\overline{a_j}=0$ then this derivative equals $x_j^{d_j-1}$ and vanishes if and only if $x_j=0$. 

 Suppose there is a $j$ such that $a_j \neq 0$ and $x_j=0$. Then for all 
 for all $k\neq j$ we have
 \[ 0=\frac{\partial F_{\overline{\lambda}}}{\partial x_k}=\overline{d_k}x_k^{d_k-1}+ \overline{a_k\lambda} \frac{\prod x_i^{a_i}}{x_k}=\overline{d_k}x_k^{d_k-1}.\] This implies that all the $x_k$ would vanish. Hence if $X_{\overline{\lambda}}$ is singular at $(x_0:\dots,x_n)$ then $x_j=0$ if $\overline{a_j}=0$ and $x_j \neq 0$ if $a_j \neq 0$.
  If there is a $j$ such that $p$ divides a non-zero $a_j$ then $X_{\overline{\lambda}}$ is quasi-smooth.  
  
  Suppose now that $p$ does not divide any of the positive $a_i$.
 
  Suppose $\overline{a_j}\neq 0$.  Consider now the derivative with respect to $x_j$:
 \[ \frac{\partial F_{\overline{\lambda}}}{\partial x_j}=\overline{d_j}x_j^{d_j-1}+ \overline{a_j\lambda} \frac{\prod x_i^{a_i}}{x_j}.\]
 This derivative vanishes if and only if 
 \[ \overline{\lambda} \prod x_i^{a_i}= -\frac{\overline{d_j}}{\overline{a_j}} x_j^{d_j}.\]
 In particular, we have
 \[ -\frac{\overline{d_j}}{\overline{a_j}} x_j^{d_j}=-\frac{\overline{d_k}}{\overline{a_k}} x_k^{d_k}\]
 for $j,k\in I$.
 
 Fix $d_j$-th roots $\alpha_j$ of $\overline{d_j}/\overline{a_j}$. Let $\zeta$ be a primitive $d'$-th root of unity. A solution of the above set of equations is of the form
 \[ x_j = \frac{\gamma^{w_j}}{\alpha_j} \zeta^{k_j w_j}\]
 for some $\gamma, k_j$.
 
 Substituting gives 
 \[ \overline{\lambda} \prod_{i\in I} \frac{\zeta^{k_iw_ia_i}}{\alpha^{a_i}}=-1,\]
which is equivalent with
\[ \overline{\lambda}^{d'} = (-1)^{d'}\prod \alpha^{a_jd'}= \frac{(-1)^{d'} \overline{d}^{d'}}{\prod_{i\in I} (\overline{a_iw_i})^{a_iw_i/g}}. \]
 \end{proof}

Let $H^i_{\dR}(U_{\lambda},\Q_q)$ denote the algebraic de Rham cohomology of $U_{\lambda}$.
\begin{Thm}[Baldassarri-Chiarellotto]
\label{deRhamcomp} Suppose $\lambda$ is chosen such that $X_{\overline{\lambda}}$ is quasi-smooth. Then
the natural map
\[ H^i_{\dR}(U_{\lambda},\Q_q) \to H^i(U_\lambda,\Q_q) \]
is an isomorphism.
\end{Thm}

\begin{proof} Consider first the case $\Ps=\Ps^n$. Then this is precisely the main theorem of \cite{BC}.

The general case can be obtained as follows: consider the standard quotient map $\pi: \Ps^n \to \Ps$  sending $x_i$ to $x_i^{w_i}$. Let $Y_\lambda$ be $\pi^{-1}(X_\lambda)$. Let $G$ be the group associated with $\pi$.  From Lemma~\ref{lemQS} it follows that $X_\lambda$ is quasi-smooth if and only if $Y_\lambda$ is smooth. Let $V_\lambda$ be the complement of $Y_\lambda$ in $\Ps^n$. Then we have an isomorphism
\[ H^i_{\dR}(V_{\lambda},\Q_q) \to H^i(U_\lambda,\Q_q). \]
There is a natural $G$-action on both groups and it is easy to see that this isomorphism is $G$-equivariant. Moreover, using \cite[Lemma 2.2.2]{Dol} we obtain that $\pi$ induces  isomorphisms $H^n_{\dR}(V_\lambda)^G\cong H^n_{\dR}(U_\lambda)$ and $H^n(V_\lambda)^G\cong H^n(U_\lambda)$, hence the natural map 
\[ H^n_{\dR}(U_{\lambda},\Q_q) \to H^n(U_\lambda,\Q_q) \]
is an isomorphism.
\end{proof}


\begin{Thm}\label{thmLefschetzHyperplane}  Let $\lambda_0 \in \Z_q$ be such that $X_{\overline{\lambda_0}}$ is quasi-smooth.
 Then the groups $H^i(U_{\lambda_0},\Q_q)$ are zero except for $i=0,n$.
\end{Thm}
\begin{proof} 
It is well-known that for a complement of a quasi-smooth hypersurface we have that $H^i_{\dR}(U,\Q_q)=0$ for $i\neq 0,n$. Applying the previous Theorem yields the result.
\end{proof}

\begin{Thm}[{Poincar\'e-Duality for $H^i(U_{\lambda},\Q_q$)}]\label{thmPoincare} Let $\overline{\lambda_0} \in \F_q$ be such that $X_{\overline{\lambda_0}}$ is quasi-smooth.
 There is a non-degenerate pairing
 \[ H^i_{\rig,c}(U_{\overline{\lambda_0}},\Q_q) \times H^{2n-i}_{\rig}(U_{\overline{\lambda_0}},\Q_q) \to H^{2n}_{\rig,c}(U_{\overline{\lambda_0}},\Q_q) \]
  respecting the Frobenius-action.
  \end{Thm}

\begin{proof}
Set $X=X_{\lambda_0}$ and $U=U_{\lambda_0}$. Consider first the case  $\Ps=\Ps^n$. Then from Lemma~\ref{lemQS} it follows that $U$ and $X$ are smooth. The main theorem of  \cite{BerPoi} asserts the existence of such pairings.

The general case can be obtained as follows: consider the standard quotient map $\pi: \Ps^n \to \Ps$  sending $x_i$ to $x_i^{w_i}$. Let $Y$ be $\pi^{-1}(X)$. Let $G$ be the group associated with $\pi$.  From Lemma~\ref{lemQS} it follows that $X$ is quasi-smooth if and only if $Y$ is smooth. Let $V$ be the complement of $Y$ in $\Ps^n$. Since Poincar\'e duality is $G$-equivariant, one obtains a pairing 
\[ H^i_{\rig}(V,\Q_q)^G \times H^{2n-i}_{\rig,c}(V,\Q_q)^G \to H^{2n}_{\rig,c}(V,\Q_q).\]

Using the isomorphism $(\Omega^k_V)^G\cong \Omega^k_U$, we obtain isomorphisms \[H^i_{\rig}(V,\Q_q)^G\cong H^i_{\rig}(U,\Q_q)\mbox{ and }H^{2n-i}_{\rig,c}(V,\Q_q)^G\cong H^{2n-i}_{\rig,c}(U,\Q_q).\] This yields the proof. 
\end{proof}

\begin{Thm}[{Lefschetz Trace Formula}] \label{thmLTF}  Let $\overline{\lambda_0} \in \F_q$ be such that $X_{\overline{\lambda_0}}$ is quasi-smooth.
Then we have that
\[ \sum_i (-1)^i \trace((q^n (\Frob^*)^{-1}) | H^i(U_{\lambda_0})) = \# U_{\overline{\lambda_0}}(\F_q).\]
\end{Thm}

\begin{proof}
Combine the Lefschetz Trace Formula for rigid cohomology with compact support \cite[Th\'eor\`eme I]{EtSt} with Poincar\'e Duality (Theorem~\ref{thmPoincare}) and Proposition~\ref{prpCompare}.
\end{proof}

 \begin{Prop} The group $H^0(U_{\lambda},\Q_q)$ is one-dimensional, and Frobenius acts trivially on $H^0(U_{\lambda},\Q_q)$.\end{Prop}
 
 \begin{proof} Straightforward.
 \end{proof}

Let $G$ be the defining equation of a quasi-smooth hypersurface $Y \subset \Ps$. Let $V:=\Ps\setminus Y$.  
Similar to the case of ordinary projective space, the algebraic de Rham cohomology of $V$ can be computed using the complex
$C^p_k = \Omega^p((k+p)Y)$. I.e., the hypercohomology group $\HH^n(\Ps,C^\bullet_k)$ equals $H^0(\Ps,C^n_k)/d H^0(\Ps,C^{n-1}_k)$ and 
\[  H^n_{\dR}(V) = \oplus_k H^0(\Ps,C^n_k)/d H^0(\Ps,C^{n-1}_k).\] 
(A proof of this equality can be obtained as follows. After fixing an embedding $\Q_q \hookrightarrow \C$ and tensoring both sides with $\C$, we obtain that it suffices to prove this result over $\C$. This is precisely the main result of  \cite{Ste}.) 

More explicitely, the vector space $H^n(V,\Q_q)$ can be identified with the quotient of the infinite-dimensional vector space spanned by
\[ \frac{H}{G^t} \Omega\]
with $\deg(H)=t\deg(G)-\sum w_i$, by the relations
\[ \frac{(t-1)H G_x-GH_x}{G^t} \Omega, \]
where the subscript $x$ means the partial derivative with respect to  a coordinate $x$ on $\Ps$. 

If $G=F$ (the polynomial whose zero-set is the Fermat hypersurface) then this formula reads as
\[ \frac{(t-1)d_i H x_i^{d_i-1} }{F^t} \Omega = \frac{H_{x_i}}{F^{t-1}}\Omega\]
in $H^n(U)$. This motivates the following definition:
\begin{Def}\label{DefRed} Let $\omega \in \Omega^n(U_0)$ be a form of the type
\[ \frac{H}{F^t} \Omega \]
with $H$ a monomial.
Let $x_i$ be a coordinate of $\Ps$ such that $x_i^{d_i-1}$ divides $H$. Then the {\em reduction  of $\omega$ with respect to $x_i$} is the form
\[ \frac{ \frac{\partial}{\partial x_i} \left(\frac{H}{x_i^{d_i-1}}\right) }{ (t-1) d_i F^{t-1} } \Omega.\]

The complete reduction $\red \omega:=H'/F^s\Omega$ of $\omega$ is the form  obtained by successively reducing with respect to the coordinates  $x_i$ of $\Ps$, such that for all $i$ the exponent of $x_i$ in $H'$ is at most $d_i-2$.
\end{Def}
Note that $\omega$ and the reduction with respect to $x_i$ of $\omega$ represent the same class in $H^n(U_0,\Q_q)$, and that the complete reduction of $\omega$ cannot be further reduced.

 \begin{Def}
 Let $P^\bullet$ be the pole order filtration on $H^n(U_\lambda)$, that is $\omega\in P^t$  if $\omega=\frac{G}{F^t_{\lambda}} \Omega$ for some $G\in \Q_q[x_0,\dots x_n]$.
 \end{Def}
 
 Let  $\mathbf{k}$ be an admissible monomial type. Recall that we can associate a differential form $\omega_\kv$ with it. By definition we have that $\omega_\kv \in P^t$, where $t$ is the relative degree of $\kv$.

 \begin{Prop} \label{prpBasis} Let $\lambda$ be such that $X_{\lambda}$ is quasi-smooth. Then the set
 \[  \left\{ \omega_\kv \colon \kv \mbox{ an admissible monomial type}\right\}\]
  is a basis for $H^n(U_{\lambda},\Q_q)$.
 \end{Prop}
\begin{proof}

The above discussion implies the statement for $\lambda=0$. 

 We start by proving  
 that for every integer $t$ the set
\[\{w_{\kv} : \kv \mbox{ an admissible monomial type of relative degree } t\}\]
is linearly independent in $P^t/P^{t-1}$.

The relations in $P^t/P^{t-1}$ are generated by (cf. the discussion before Definition~\ref{DefRed}) 
\[ \frac{x_i^{d_i-1}\prod x_j^{k_j}}{F^t_{\lambda}}\Omega= \frac{-\lambda a_i}{d_i} \frac{\prod x_j^{k_j+a_j}}{x_iF^t_{\lambda}}\Omega.\]
Suppose $i$ is chosen such that $a_i\neq 0$. Let 
\[ \sigma_i(G):=\frac{-d_i x_i^{d_i}}{\lambda a_i\prod_jx_j^{a_j}}G.\]
If $G$ is a monomial of degree $td-\sum w_i$ such that all the exponents of the $x_j$ are at least $a_j-\delta_{i,j}$, then
\[ \frac{G}{F^t_{\lambda}}\Omega \equiv \frac{\sigma_i(G)}{F^t_{\lambda}} \Omega \bmod P^{t-1}.\]
Note that $\sigma_i$ is defined if the exponent of $x_j$ is at least $a_j-\delta_{i,j}d_i$, but $\sigma_i$ corresponds to a relation in $P^t$ only if the exponent $x_j$ is at least $a_j-\delta_{i,j}$.
Similarly, if the exponent of $x_i$ in $G$ is at least $d_i-1$, then 
\[ \frac{G}{F^t_{\lambda}}\Omega \equiv \frac{\sigma_i^{-1}(G)}{F^t_{\lambda}} \Omega \bmod P^{t-1}.\]
 
Take a non-trivial expression $\sum b_{\kv} \omega_{\kv}$ that is zero modulo $P^{t-1}$.
Since the $\sigma_i$ generate the relations, and the $\sigma_i$ map monomials to monomials, there exists two distinct admissible monomial types $\kv,\mv$ of relative degree $t$ and a sequence of $\sigma_i$ and $\sigma_j^{-1}$ such that\[ \tau(\omega_{\kv}):=\sigma_{i_s}^{\epsilon_s}\dots \sigma_{i_1}^{\epsilon_1}(\omega_{\kv})= c_0\omega_{\mv}, \]
 with
  \begin{itemize}
\item  $c_0 \in \Q_q$, 
\item  $\epsilon_j\in \{ \pm 1\}$,
\item $a_{i_j}\neq 0$ for all $j$, 
\item  for all $j$ such that $\epsilon_j=1$, we have that for all $k$ the exponent of $x_k$ in $\sigma_{i_{j-1}}^{\epsilon_{j-1}}\dots \sigma_{i_1}^{\epsilon_1}(\omega_{\kv})$ is at least  $a_k-\delta_{i_j,k}$ and
\item  for all $j$ such that $\epsilon_j=-1$, we have that  the exponent of $x_{i_j}$ in $\sigma_{i_{j-1}}^{\epsilon_{j-1}}\dots \sigma_{i_1}^{\epsilon_1}(\omega_{\kv})$ is at least  $d_{i_j}-1$.
\end{itemize}
 We will prove below that given such a $\tau$, we can always shorten the length of this expression by $2$, and that this expression cannot consist of one $\sigma_i$. Hence the only possibility for $\tau$ is to be the identity and $k_i=m_i$ for all $i$, a contradiction. 

We claim that $\epsilon_1=1$ and $\epsilon_s=-1$. If $\epsilon_1$ were $-1$, then in order to apply $\sigma_{i_1}$ we would need that  the exponent $x_{i_1}$ in $\omega_{\kv}$ is at least $d_{i_1}-1$, contradicting that $\omega_{\kv}$ is associated with an admissible monomial type. Similarly, if $\epsilon_s=1$ we obtain that the exponent of $x_{i_s}$ in $\omega_{\mv}$ is at least $d_{i_s}-1$, contradicting that $\omega_{\mv}$ is associated with an admissible monomial type.
 
 %


Let $j$ be the smallest integer such that $\epsilon_j=-1$. This implies that the exponent of $x_{i_j}$ in $ \sigma_{i_{j-1}}\dots \sigma_{i_1}(\prod x_i^{k_i})$ is at least $d_{i_j}-1$, hence at least for one of the $j'<j$ we have $i_j=i_{j'}$. Let $j'$ be the largest integer smaller than $j$ such that $i_j=i_{j'}$.
 
Note that the $\sigma_i$ commute as operators on $\Q_q(x_0,\dots,x_n)$. Hence, if we consider the $\sigma_i$ as operators on $\Q_q(x_0,\dots x_n)$ then we have the identities
 \begin{eqnarray*} \sigma_{i_j}^{-1}\sigma_{i_{j-1}} \dots \sigma_{i_{j'}} \sigma_{i_{j'-1}} \dots \sigma_{i_1}\left(\prod x_i^{k_i}\right)&=& \sigma_{i_j}^{-1}\sigma_{i_{j'}}\sigma_{i_{j-1}} \dots \sigma_{i_{j'+1}}\sigma_{i_{j'-1}} \dots \sigma_{i_1}\left(\prod x_i^{k_i}\right)\\&=&\sigma_{i_{j-1}} \dots \sigma_{i_{j'+1}}\sigma_{i_{j'-1}} \dots \sigma_{i_1}\left(\prod x_i^{k_i}\right).\end{eqnarray*}
We need to show that the latter expression corresponds to a series of relations in $P^t/P^{t-1}$, i.e., 
we need to show that for each $j''$ such that $j'<j''<j$ we have that if 
\[ \sigma_{i_{j''}} \dots \sigma_{i_{j'+1}}\sigma_{i_{j'-1}} \dots \sigma_{i_1}
\left(\prod x_i^{k_i}\right) = c \prod x_r^{e_r}\]
with $c\in \Q_q$
then $e_r\geq  a_r-\delta_{r,i_{j''+1}}$ for all $r$.

Suppose that $r\neq i_j$. Since  
\[  \sigma_{i_{j''}}  \dots \sigma_{i_1}             \left(\prod x_i^{k_i}\right)= c' \prod x_r^{e'_r} \]
with $c'\in \Q_q$ and $e'_r\geq a_r-\delta_{r,i_{j''+1}}$ and  $\sigma_{i_j}$ lowers the exponent of $x_r$ by $a_r$ we obtain that  $e'_r=e_r-a_r$, whence  $e_r \geq  a_r-\delta_{r,i_{j''+1}}$.
 
Suppose that $r=i_j$. Since 
 \[ \sigma_{i_{j-1}} \dots \sigma_{i_1} \left(\prod x_i^{k_i}\right)=c''  \prod x_r^{e''_r} \]
with $c''\in \Q_q$ and $e''_r\geq d_r-1$, it follows that 
 \[ \sigma_{i_{j-1}} \dots \sigma_{i_{j'+1}}\sigma_{i_{j'-1}} \dots \sigma_{i_1} \left(\prod x_i^{k_i}\right)=c'''  \prod x_r^{e'''_r} \]
with $c'''\in \Q_q$ and $e'''_r\geq 0$. Since the $\sigma_{i_k}$ for $j''<k<j'$ lower the exponent of $x_r$ by $a_r$ we obtain that
  $e_r=e'''_r+(j-j'')a_r\geq a_r$. 

We need to show that \[  \left\{ \omega_\kv \colon \kv \mbox{ an admissible monomial type}\right\}\]
 spans $H^n(U_{\lambda},\Q_q)$. If $\lambda=0$ then this follows from the discussion before this proposition. If $\lambda\neq 0$ and all the weights equal 1 then  \cite[Theorem 1.10]{Katz} shows that $\dim H^{n-1}(X_\lambda,\Q_q)$ is independent of $\lambda$. Using tzhe Gysin exact sequence we obtain that $\dim H^n(U_\lambda,\Q_q)$ is independent of $\lambda$. The general case follows from this case by applying the standard quotient map and \cite[Lemma 2.2.2]{Dol}.
\end{proof}

\section{Deformation Theory}\label{SecDeform}
Assume for the moment that $\Ps=\Ps^n$. Following N. Katz, consider the  commutative diagram
\[\xymatrix {H^n(U_{\lambda^q}) \ar[r]^{\Frob_{q,\lambda}^*} \ar[d]_{A(\lambda^q)}& H^n(U_\lambda) \ar[d]_{A(\lambda)} \\ H^n(U_0) \ar[r]^{\Frob_{q,0}^*} & H^n(U_0),}\]
where  $\Frob_{q,\lambda}$ is the Frobenius acting on the complete family. Since it maps the fiber over $0$ to the fiber over $0$ this map can be restricted to $U_0$.  Katz studied the differential equation  associated to $A(\lambda)$. He remarked in a note that $A(\lambda)$ is actually the solution of the Picard-Fuchs equation. 

We first give a way of computing a map $B(\lambda)$ such that $\Frob_{q,0}^*B(\lambda^q)=B(\lambda)\Frob_{q,\lambda}$ on a small neighborhood of $0$. This matrix $B(\lambda)$ is enough to deduce $\Frob_{q,\lambda}^*$ from $\Frob_{q,0}^*$. 

Fix a basis \[ \frac{G_i}{F_\lambda^{t}} \Omega\]
for $H^n(U_\lambda)$ and  write
\begin{equation}\label{maineqn} \frac{G_i}{F_\lambda^{t}} \Omega= \sum_{j=0}^\infty \left( \begin{array}{c} j+t-1 \\ j \end{array} \right) \frac{G_i (F-F_\lambda)^j}{F^{j+t}} \Omega.\end{equation}
Since $F-F_\lambda$ is the product of $\lambda$ with a polynomial with integral coefficients, we have that the above power series in the $x_i$ converges on a small disc. By choosing $\lambda$ sufficiently small, we obtain an overconvergent power series in the $x_i$, hence $\frac{G_i}{F_\lambda^t}$ defines an element of  $H^n(U_0)$. Let $B(\lambda) :  H^n(U_\lambda) \rightarrow H^n(U_0)$ be the analytic continuation of the operator mapping
\[ \frac{G_i}{F_\lambda^{t}} \Omega \] to the complete reduction of (\ref{maineqn}) in $H^n(U_0)$.

In this way we obtain a local expansion of the matrix $B(\lambda)$ around $\lambda=0$. In the following section we will make this more explicit.

\begin{Prop}[Katz] \label{prpKatz} We have $B(\lambda)\Frob_{q,\lambda}^*=\Frob_{q,0}^*B(\lambda^q)$ and $B(\lambda)=A(\lambda)$.
\end{Prop}
\begin{proof} The case $\Ps=\Ps^n$  is a combination of \cite[Lemma 2.10]{Katz}, \cite[Lemma 2.13]{Katz} and \cite[Theorem 2.14]{Katz}. The general case is a formal consequence of the special case by Lemma~\ref{lemQS}, Proposition~\ref{prpKatz} and the definition of $H^n(U_{\lambda},\Q_q)$ in terms of   
 the standard quotient map $\pi:\Ps^n\to \Ps$.\end{proof}


\begin{Rem} Proposition~\ref{prpKatz} is particularly interesting in the case when we specialize to $\lambda=\lambda_0$ where $\lambda_0$ is the Teichm\"uller lift of some element $\overline{\lambda_0}$. Then $\lambda_0^q=\lambda_0$, hence $\Frob^*_{q,\lambda_0}$ is a lift of Frobenius  on $H^n(U_{\lambda_0},\Q_q)$. Using Theorem~\ref{thmLTF} and Theorem~\ref{thmLefschetzHyperplane} we obtain that
\[ Z(U_{\overline{\lambda_0}},t)=\lim_{\lambda \to \lambda_0} \frac{\left(\det (I-tq^n A(\lambda^q)^{-1} (\Frob^*_{q,0})^{-1} A(\lambda)) \right)^{(-1)^{n+1} }}{1-q^nt}.\]
\end{Rem}

\section{Actual computation of the deformation matrix}\label{SecComp}
In order to compute the matrix $A(\lambda)$ we need to reduce the right hand side of (\ref{maineqn}) in $H^n(U_0)$. We start with a very useful lemma.

\begin{Lem}\label{redformLem} Fix non-negative integers $b_i$ such that $\sum b_iw_i+w=td$ for some integer $t$.
The complete reduction of $\omega:=\frac{\prod x_i^{b_i}}{F^t}\Omega$ equals
\[ \frac{\prod_i ((c_i+1)w_i/d)_{q_i}} {(s)_{t-s}} \frac{\prod x_i^{c_i}}{F^{s}}\Omega,\]
where $0\leq c_i< d_i$  and $q_i,s$ are integers such that $b_i= q_i d_i+c_i$, and $sd=\sum c_iw_i+w$, i.e., $t-s=\sum q_i$, provided that $c_i \neq d_i - 1$ for all $i$. If for one of the $i$ we have $c_i=d_i-1$ then $\omega$ reduces to zero in $H(U_0,\Q_q)$. 
\end{Lem}
\begin{proof}

The reduction with respect to $x_0$ of
\[ \frac{x_0^{b_0}\prod_{i=1}^n x_i^{b_i}}{F^t}\Omega\]
(cf. Definition~\ref{DefRed}) equals
\[ \frac{ x_0^{b_0-d_0} ((b_0+1)-d_0) \prod x_i^{b_i}}{(t-1)d_0F^{t-1}}\Omega=\frac{ x_0^{b_0-d_0} ((b_0+1)w_0-d) \prod x_i^{b_i}}{(t-1)dF^{t-1}}\Omega\]
(provided  $b_0\geq d_0$). After reducing $q_i$ times  with respect to $x_i$ for $i=0,\dots,n$, we obtain that  $\omega$ reduces to  
\[ \frac{(s-1)! \prod_i \left( \prod_{j=0}^{q_i-1}((c_i+1)w_i+jd)    \right)\prod x_i^{c_i}}{(t-1)! d^{t-s} F^s} \Omega.\]
This in turn equals
\[ \tau:= \frac{(s-1)! d^{\sum q_i} \prod ((c_i+1)w_i/d)_{q_i} \prod x_i^{c_i}} {(t-1)! d^{t-s} F^s}\Omega.\]
If none of the $c_i$ equals $d_i-1$ then this is a complete reduction. Using $\sum q_i=t-s$ the first formula follows.

If $c_i=d_i-1$ then we can write $\tau$ as $(F_{x_i}G/F^s) \Omega$, where $G$ does not contain the variable $x_i$.  The reduction of this form is a constant times 
\[ \frac{G_{x_i}}{F^{s-1}} \Omega. \]
Since $G_{x_i}=0$, this reduction is zero.
\end{proof}

Fix an admissible monomial type $\mathbf{k}=(w_0(k_0+1),\dots,w_n(k_n+1))\in (\Z/d\Z)^{n+1}$ of relative degree $t$. We want to calculate the reduction of
\[ \frac{ \prod x_i^{k_i}}{(F+\lambda \prod x_i^{a_i})^t} \Omega\]
in $H^n(U_0)$. In order to find a power series expression, we assume that $\lambda$ is sufficiently small, then by (\ref{maineqn}) this form equals
\begin{equation}\label{eqndet} \frac{\prod x_i^{k_i}}{F^t} \frac{1}{\left(1-\frac{(-\lambda) \prod x_i^{a_i}}{F}\right)^t} \Omega= \sum_j \left( \begin{array}{c} t+j-1 \\ j \end{array}\right) \frac{\prod x_i^{k_i+a_i j}}{F^{t+j}} (-\lambda)^j \Omega.\end{equation}

Note that at most $d$ distinct monomials occur in the reduction of the form. 

\begin{Def}
Let  $r,s$ be non-negative integers, let $\alpha_i \in \Q_q$, for $i \in \{1,2,\dots, r\}$, let $\beta_j \in \Q_q \setminus \Z_{<0}$ for $j\in \{1,2,\dots, s\}$. We define the \emph{(generalized) hypergeometric function}
\[{}_{r}F_{s}\left( \begin{array}{c} \alpha_{1} \;\alpha_2\; \dots \; \alpha_r  \\ \beta_1 \;\beta_2 \; \dots \; \beta_s  \end{array}; z    \right)\] to be
\[ \sum_{k=0}^\infty b_j z^j,\]
with $b_0=1$, and
\[ \frac{b_{j+1}}{b_j} = \frac{(j+\alpha_1)\dots (j+\alpha_r)}{(j+\beta_1)\dots (j+\beta_s) (j+1)},\] for all positive integers $j$.
\end{Def}

Let 
$d'_i$ be  the order  of $a_i \bmod d_i$ in $ \Z/d_i\Z$. Let $d'$ be the least common multiple of all the $d_i'$.
 Set $b_i = a_i d'/d_i$. In the following Proposition and its proof we identify elements in $a \in \Z/m\Z$ with their representative $\tilde{a}\in \Z$ such that $0\leq \tilde{a} \leq m-1$.
 
\begin{Prop}\label{PrpCoeff}Let $\kv$ be an admissible monomial type. Let $t$ be the relative degree of $\kv$. Write  $A(\lambda)\omega_\kv =\sum c_\mv(\lambda) \omega_\mv$, where  the sum is taken over all admissible monomial types. Then $c_\mv(\lambda)$  is non-zero only if   
there is a $j_0\in \Z$  with $0\leq j_0\leq d'-1$ and such that $\mathbf{m}-\mathbf{k}=\overline{j_0} \mathbf{a}$. If this is the case then $c_{\mv}(\lambda) \omega_\mv$ equals $c'_{\mv}(-\lambda)\red \frac{\prod x_i^{a_i j_0+k_i} }{F^{t+j_0}} \Omega$,
where
\[c'_{\mv}(\lambda):= \left( \begin{matrix} t+j_0-1 \\ j_0 \end{matrix}\right) \lambda^{j_0}  {}_{d'}F_{d'-1}\left( \begin{array}{c} \alpha_{i,s} \\ \frac{j_0+1}{d'} \; \frac{j_0+2}{d'} \; \dots \; \widehat{1} \; \dots \;\frac{j_0+d'}{d'} \end{array}; \prod_{i\colon a_i\neq 0} \left(\frac{  a_i}{d_i}\right)^{b_i} \lambda^{d'}    \right) \]
times $\red \frac{\prod x_i^{a_i j_0+k_i} }{F^{t+j_0}} \Omega$,
with
\[ \alpha_{i,s} = \frac{(s-1)d_i+1+a_i j_0+k_i}{a_i d'}, s=1,\dots,b_i; i=0,\dots n.\]
\end{Prop}
Note that this Proposition almost gives  a complete reduction of the differential form $\Frob_{q,\lambda}(\omega_\kv)$, in the sense that $c_{\kv}(\lambda)$ is described as the product of a hypergeometric function and the reduction of a rational function in the $x_i$ multiplied by $\Omega$. The latter form can be easily reduced using Lemma~\ref{redformLem}.
\begin{proof}
It suffices to   compute explicitly a complete reduction of
\[ \omega:=\frac{\prod x_i^{k_i}}{F_\lambda}\Omega \]
in $H^n(U_0)$. We can write $\omega$ as
\[  \sum_j \left( \begin{array}{c} t+j-1 \\ j \end{array}\right) \frac{\prod x_i^{k_i+a_i j}}{F^{t+j}} (-\lambda)^j \Omega.\] 
Set $c_{t,j}:=\left( \begin{array}{c} k+j-1 \\ j \end{array}\right)$. Since each reduction  step decreases the exponent of $x_i$ by $d_i$, we split this sum as follows: 
write
\[  \omega= \sum_{j_0=0}^{d'-1} \sum_j c_{t,j_0+d' j} \frac{\prod x_i^{k_i+a_i (j_0+d'j)}}{F^{t+j_0+d'j}} (-\lambda)^{j_0+d'j}\Omega.\] 

For $0\leq j_0 \leq d'-1$ set 
\[ \omega_{j_0} := \sum_j c_{t,j_0+d' j} \frac{\prod x_i^{k_i+a_i (j_0+d'j)}}{F^{t+j_0+d'j}} (-\lambda)^{j_0+d'j}\Omega.\]

From Lemma~\ref{redformLem} it follows that if  for some $i$ we have that $k_i+a_i(j_0+d'j)\equiv -1 \bmod d_i$ then $\omega_{j_0}$ reduces to zero. Otherwise, we claim that the reduction of $\omega_j$ is a generalized hypergeometric function. In order to prove this and to calculate the parameters, we need to show that 
\begin{equation}\label{eqnquo} \frac{ c_{t,j_0+d'j+d'} \red \frac{ \prod x_i^{k_i+a_i (j_0+d'j) +a_i d'}}{F^{t+j_0+d'j+d'}} \Omega}{ c_{t,j_0+d' j} \red \frac{\prod x_i^{k_i+a_i (j_0+d'j)}}{F^{t+j_0+d'j}}\Omega}
\end{equation}
is a rational function in $j$. If we reduce with respect to $x_i$ then the exponent of $x_i$ is lowered by $d_i$. So if we reduce the numerator $b_i=a_id'/d_i$ times with respect to $x_i$, then the exponent of $x_i$ in the numerator and denominator coincide. Now
\[ \red \frac{ \prod x_i^{k_i+a_i (j_0+d'j) +a_i d'}}{F^{t+j_0+d'j+d'}} \Omega\] equals \[\frac{\prod_i  \prod_{s=1}^{b_i} (k_i+a_i(j_0+d'j)+(s-1) d_i+1)}{(t+j_0+d'j)_{\sum b_i} \prod_i d_i^{b_i}}\red \frac{\prod x_i^{k_i+a_i (j_0+d'j)}}{F^{t+j_0+d'j}} \Omega
\]
and
\[ \frac{ c_{t,j_0+d'j+d'}}{ c_{t,j_0+d' j} } =  \frac{(t+j_0+d'j)_{d'}}{(j_0+d'j+1)_{d'}}.\]

Putting this together we obtain that (\ref{eqnquo})  equals:
\[\frac{\prod_i  \prod_{s=1}^{b_i} (k_i+a_i(j_0+d'j)+(s-1) d_i+1)}{(j_0+d'j+1)_{d'} \prod_i d_i^{b_i}}.\]

This equals
\[  \frac{\prod_{i:a_i\neq 0}(a_id')^{b_i} }{(d')^{d'}\prod_{i: a_i\neq 0} d_i^{b_i} }  \frac{\prod_i\prod_{s=1}^{b_i} (j + \frac{(s-1)d_i+1+a_i j_0+k_i}{a_i d'} )} {\prod_{s=1}^{d'} (j+\frac{j_0+s}{d'})}.\]
Since $\sum b_i=d'$ the first factor simplifies to
\[ \prod_{i:a_i\neq 0}\left(\frac{a_i}{d_i}\right)^{b_i}. \]
From the second factor we can read off the hypergeometric parameters.

Since the first summand of $\omega_{j_0}$ equals
\[ (-\lambda)^{j_0} c_{t,j_0}\frac{\prod x_i^{a_i j_0+k_i} }{F^{t+j_0}} \Omega,\]
by collecting everything together, we obtain that \[
\frac{\red \omega_j}{c_{t,j_0}(-\lambda)^{j_0}}\] equals
\[     {}_{d'}F_{d'-1}\left( \begin{array}{c} \alpha_{i,s} \\ \frac{j_0+1}{d'} \; \frac{j_0+2}{d'} \; \dots \; \widehat{1} \; \dots \;\frac{j_0+d'}{d'} \end{array}; ( - \lambda )^{d'}\prod_{i:a_i\neq0}    \left(\frac{a_i}{d_i}\right)^{b_i} \right) \red \frac{\prod x_i^{a_i j_0+k_i} }{F^{t+j_0}} \Omega,\]
as desired.
\end{proof}

\begin{Exa}  Consider the family $X^3+Y^3+Z^3+\lambda XYZ$. Then we obtain the following matrix $A(\lambda)$  (with respect to the basis $\{\omega_{(1,1,1)},\omega_{(2,2,2)}\}$)
\[ \left(\begin{array}{cc} 
 \;_{2}F_1\left( \begin{array}{c} \frac{1}{3} \; \frac{1}{3} \\ \frac{2}{3} \end{array} ; \frac{-\lambda^3}{27}
\right)  & 
\frac{\lambda^2}{54} \;_{2}F_1\left( \begin{array}{c} \frac{4}{3} \; \frac{4}{3} \\ \frac{5}{3} \end{array} ; \frac{-\lambda^3}{27}
\right)  \\
-\lambda \;_{2}F_1\left( \begin{array}{c} \frac{2}{3} \; \frac{2}{3} \\ \frac{4}{3} \end{array} ; \frac{-\lambda^3}{27}
\right)  & 
 \;_{2}F_1\left( \begin{array}{c} \frac{2}{3} \; \frac{2}{3} \\ \frac{1}{3} \end{array} ; \frac{-\lambda^3}{27}
\right)  

\\ \end{array}\right). \]
\end{Exa}

\begin{Exa}
Another famous example is $C_\lambda:X^4+Y^4+Z^4+\lambda X^2Y^2$.

Note that $d'=2$, $d_1'=d_2'=2, d'_3=0, b_1=b_2=1, b_3=0$. 

One easily obtains the following results
\begin{eqnarray*} 
A(\lambda)\omega_{(2,1,1)}&=& _1F_0\left( \begin{matrix} \frac{1}{4} \\ - \end{matrix} ; \frac{\lambda^2}{16} \right)\omega_{(2,1,1)},\\
A(\lambda)\omega_{(1,2,1)}&=& _1F_0\left( \begin{matrix}\frac{1}{4} \\ - \end{matrix} ; \frac{\lambda^2}{16} \right)\omega_{(1,2,1)},\\
A(\lambda)\omega_{(2,3,3)}&=& _1F_0\left( \begin{matrix} \frac{3}{4}\\ - \end{matrix} ; \frac{\lambda^2}{16} \right)\omega_{(2,3,3)},\\
A(\lambda)\omega_{(3,2,3)}&=& _1F_0\left( \begin{matrix} \frac{3}{4} \\ - \end{matrix} ; \frac{\lambda^2}{16} \right)\omega_{(3,2,3)}.
\end{eqnarray*}

$A(\lambda)$ acts as follows on the basis $\{\omega_{(1,1,2)},\omega_{(3,3,2)}\}$
\[ \left(\begin{array}{cc} 
 \;_{2}F_1\left( \begin{array}{c} \frac{1}{4} \; \frac{1}{4} \\ \frac{1}{2} \end{array} ; \frac{\lambda^2}{16}
\right)  & 
\frac{\lambda^2}{16} \;_{2}F_1\left( \begin{array}{c} \frac{5}{4} \; \frac{5}{4} \\ \frac{3}{2} \end{array} ; \frac{\lambda^2}{16}
\right)  \\
-\lambda \;_{2}F_1\left( \begin{array}{c} \frac{3}{4} \; \frac{3}{4} \\ \frac{3}{2} \end{array} ; \frac{\lambda^2}{16}
\right)  & 
 \;_{2}F_1\left( \begin{array}{c} \frac{3}{4} \; \frac{3}{4} \\ \frac{1}{2} \end{array} ; \frac{\lambda^2}{16}
\right)  
\\ \end{array}\right). \]

It is classically known that the Jacobian of $C_\lambda$ is isogenous to the product of two elliptic curves with $j$-invariant $1728$  and one elliptic curve $E_\lambda$ whose $j$-invariant depends properly on $\lambda$. This factor can also be obtained  from the above information:

When we restrict $A(\lambda)$ to the subspace spanned by $\omega_{(1,1,2),(3,2,2)}$, we find the same operator as the operator $A'(\lambda)$ associated with the family $E_\lambda:X^4+Y^4+Z^2+\lambda X^2Y^2$  considered in $\Ps(1,1,2)$. One easily shows that this is a family of elliptic curves, with $j$-invariant depending on $\lambda$.
The curve $E_0$ has an automorphism of order 4 with fixed points, hence $j(E_0)=1728$.

In the next section we prove that if $q\equiv 1 \bmod 4$ then all the $\omega_\kv$ are eigenvectors for $\Frob^*_{q,0}$, let $c_{\kv,q}$ be the corresponding eigenvalue.  Then for $\kv=(2,1,1)$ we have that 
\[ \Frob^*_{q,\lambda} \omega_\kv=A(\lambda)^{-1} \Frob^*_{q,0} A(\lambda^q)\omega_\kv= \frac{_1F_0\left( \begin{matrix} \frac{1}{4} \\ - \end{matrix} ; \frac{\lambda^{2q}}{16} \right)}{_1F_0\left( \begin{matrix}\frac{1}{4} \\ - \end{matrix} ; \frac{\lambda^2}{16} \right)} c_{\kv,q}\omega_{\kv}.\]
One easily shows that the factor in front of $c_{\kv,q}$ is a 4-th root of unity, which implies that we have twisted the Frobenius action  on $\omega_\kv$  by a quartic character. Something similar happens when $\kv\in \{ (1,2,1),(2,3,3),(3,2,3)\}$. This implies that on a 4-dimensional subspace $V_\lambda$ of $H^1(X_\lambda,\Q_q)$ the Frobenius-action is a quartic twist of the Frobenius-action on $V_0\subset H^1(X_0,\Q_q)$. 
The curve $X_0$ has the automorphism $[X,Y,Z]\mapsto [Z,X,Y]$. From this we obtain that the action of Frobenius on $V_0$ is isomorphic to two copies of the Frobenius action on $E_0$.
\end{Exa}

\begin{Exa}\label{exaQuiDef}
Consider now the quintic threefold $X_0^5+X_1^5+X_2^5+X_3^5+X_4^5+\lambda X_0X_1X_2X_3X_4$. This family is studied for example by 
 Candelas, de la Ossa, and Rodriguez-Villegas in \cite{Cand2}. We discuss another aspect of this family in Example~\ref{exaQuiMon}.

One can distinguish between the following five types of subspaces:

We start with $V_1=\spa  \{\omega_{(1,1,1,1,1)},\omega_{(2,2,2,2,2)},\omega_{(3,3,3,3,3)},\omega_{(4,4,4,4,4)}\}$.
The corresponding matrix is

\begin{SMALL}

\[ \left(\begin{array}{cccc}

 \;_{4}F_3\left(\begin{array}{c} (\frac{1}{5})^4  \\ \frac{2}{5} \;\frac{3}{5} \;\frac{4}{5} \end{array}; \frac{-\lambda^5}{5^5}
\right)  & 
\frac{\lambda^4}{2^3\cdot 3 \cdot 5^5}  \;_{4}F_3\left(\begin{array}{c} (\frac{6}{5})^4 \\ \frac{7}{5} \;\frac{8}{5} \;\frac{9}{5} \end{array}; \frac{-\lambda^5}{5^5}
\right)  & 
\frac{-\lambda^3}{2^2\cdot 3 \cdot 5^5}  \;_{4}F_3\left(\begin{array}{c} (\frac{6}{5})^4 \\ \frac{4}{5} \;\frac{7}{5} \;\frac{8}{5} \end{array}; \frac{-\lambda^5}{5^5}
\right)  & 
\frac{\lambda^2}{2^2\cdot 3 \cdot 5^5}  \;_{4}F_3\left(\begin{array}{c} (\frac{6}{5})^4 \\ \frac{3}{5} \;\frac{4}{5} \;\frac{7}{5} \end{array}; \frac{-\lambda^5}{5^5}
\right)  \\

{-\lambda} \;_{4}F_3\left(\begin{array}{c} (\frac{2}{5})^4  \\ \frac{3}{5} \;\frac{4}{5} \;\frac{6}{5} \end{array}; \frac{-\lambda^5}{5^5}
\right)  & 
 \;_{4}F_3\left(\begin{array}{c} (\frac{2}{5})^4 \\ \frac{1}{5} \;\frac{3}{5} \;\frac{4}{5} \end{array}; \frac{-\lambda^5}{5^5}
\right)  & 
\frac{2\lambda^4}{3\cdot  5^5}  \;_{4}F_3\left(\begin{array}{c} (\frac{7}{5})^4 \\ \frac{6}{5} \;\frac{8}{5} \;\frac{9}{5} \end{array}; \frac{-\lambda^5}{5^5}
\right)  & 
\frac{-8\lambda^3}{3^3 \cdot 5^5 }  \;_{4}F_3\left(\begin{array}{c} (\frac{7}{5})^4 \\ \frac{4}{5} \;\frac{6}{5} \;\frac{8}{5} \end{array}; \frac{-\lambda^5}{5^5}
\right)  \\

{\lambda^2} \;_{4}F_3\left(\begin{array}{c} (\frac{3}{5})^4  \\ \frac{4}{5} \;\frac{6}{5} \;\frac{7}{5} \end{array}; \frac{-\lambda^5}{5^5}
\right)  & 
{-2\lambda}  \;_{4}F_3\left(\begin{array}{c} (\frac{3}{5})^4 \\ \frac{2}{5} \;\frac{4}{5} \;\frac{6}{5} \end{array}; \frac{-\lambda^5}{5^5}
\right)  & 
  \;_{4}F_3\left(\begin{array}{c} (\frac{3}{5})^4 \\ \frac{1}{5} \;\frac{2}{5} \;\frac{4}{5} \end{array}; \frac{-\lambda^5}{5^5}
\right)  & 
\frac{27\lambda^4}{2^3 \cdot 5^5}  \;_{4}F_3\left(\begin{array}{c} (\frac{8}{5})^4 \\ \frac{6}{5} \;\frac{7}{5} \;\frac{9}{5} \end{array}; \frac{-\lambda^5}{5^5}
\right)  \\

{-\lambda^3} \;_{4}F_3\left(\begin{array}{c} (\frac{4}{5})^4  \\ \frac{6}{5} \;\frac{7}{5} \;\frac{8}{5} \end{array}; \frac{-\lambda^5}{5^5}
\right)  & 
{3\lambda^2}  \;_{4}F_3\left(\begin{array}{c} (\frac{4}{5})^4 \\ \frac{3}{5} \;\frac{6}{5} \;\frac{7}{5} \end{array}; \frac{-\lambda^5}{5^5}
\right)  & 
{-3\lambda}  \;_{4}F_3\left(\begin{array}{c} (\frac{4}{5})^4 \\ \frac{2}{5} \;\frac{3}{5} \;\frac{6}{5} \end{array}; \frac{-\lambda^5}{5^5}
\right)  & 
  \;_{4}F_3\left(\begin{array}{c} (\frac{4}{5})^4 \\ \frac{1}{5} \;\frac{2}{5} \;\frac{3}{5} \end{array}; \frac{-\lambda^5}{5^5}
\right)  \\

\end{array} \right)
\]

\end{SMALL}

Here we used the shorthand 
\[ _{4}F_3\left(\begin{array}{c} a^4 \\ b \; c \; d \end{array}; \mu \right) 
\mbox{ for }
 _{4}F_3\left(\begin{array}{c} a \; a \;a \;a \\ b \; c \; d \end{array}; \mu \right). \]
The other four spaces are less interesting: on
$V_2=\spa \{\omega_{(1,1,1,3,4)}, \omega_{(4,4,4,1,2)}\}$ we have that $A(\lambda)$ acts as  
\[ \left(\begin{array}{cc} 
 \;_{2}F_1\left( \begin{array}{c} \frac{1}{5} \; \frac{1}{5} \\ \frac{2}{5} \end{array} ; \frac{-\lambda^5}{3125}
\right)  & 
\frac{\lambda^2}{500} \;_{2}F_1\left( \begin{array}{c} \frac{6}{5} \; \frac{6}{5} \\ \frac{7}{5} \end{array} ; \frac{-\lambda^5}{3125}
\right)  \\
-{4\lambda^3} \;_{2}F_1\left( \begin{array}{c} \frac{4}{5} \; \frac{4}{5} \\ \frac{8}{5} \end{array} ; \frac{-\lambda^5}{3125}
\right)  & 
 \;_{2}F_1\left( \begin{array}{c} \frac{4}{5} \; \frac{4}{5} \\ \frac{3}{5} \end{array} ; \frac{-\lambda^5}{3125}
\right)  
\\ \end{array}\right). \]

On $V_3=\spa \{\omega_{(2,2,2,1,3)}, \omega_{(3,3,3,2,4)}\}$ we have that $A(\lambda)$ acts as
\[ \left(\begin{array}{cc} 
 \;_{2}F_1\left( \begin{array}{c} \frac{2}{5} \; \frac{2}{5} \\ \frac{4}{5} \end{array} ; \frac{-\lambda^5}{3125}
\right)  & 
\frac{\lambda^4}{6250} \;_{2}F_1\left( \begin{array}{c} \frac{7}{5} \; \frac{7}{5} \\ \frac{9}{5} \end{array} ; \frac{-\lambda^5}{3125}
\right)  \\
-{2\lambda  } \;_{2}F_1\left( \begin{array}{c} \frac{3}{5} \; \frac{3}{5} \\ \frac{6}{5} \end{array} ; \frac{-\lambda^5}{3125}
\right)  & 
 \;_{2}F_1\left( \begin{array}{c} \frac{3}{5} \; \frac{3}{5} \\ \frac{1}{5} \end{array} ; \frac{-\lambda^5}{3125}
\right)  
\\ \end{array}\right). \]

On $V_4=\spa \{\omega_{(1,1,2,2,4)}, \omega_{(3,3,4,4,1)}\}$ we have that $A(\lambda)$ acts as
\[ \left(\begin{array}{cc} 
 \;_{2}F_1\left( \begin{array}{c} \frac{1}{5} \; \frac{2}{5} \\ \frac{3}{5} \end{array} ; \frac{-\lambda^5}{3125}
\right)  & 
\frac{-\lambda^3}{1875} \;_{2}F_1\left( \begin{array}{c} \frac{6}{5} \; \frac{7}{5} \\ \frac{8}{5} \end{array} ; \frac{-\lambda^5}{3125}
\right)  \\
 \frac{\lambda^2}{5} \;_{2}F_1\left( \begin{array}{c} \frac{3}{5} \; \frac{4}{5} \\ \frac{7}{5} \end{array} ; \frac{-\lambda^5}{3125}
\right)  & 
 \;_{2}F_1\left( \begin{array}{c} \frac{3}{5} \; \frac{4}{5} \\ \frac{2}{5} \end{array} ; \frac{-\lambda^5}{3125}
\right)  
\\ \end{array}\right). \]

On $V_5=\spa \{\omega_{(3,3,1,1,2)}, \omega_{(4,4,2,2,3)}\}$ we have that $A(\lambda)$ acts as
\[ \left(\begin{array}{cc} 
 \;_{2}F_1\left( \begin{array}{c} \frac{1}{5} \; \frac{3}{5} \\ \frac{4}{5} \end{array} ; \frac{-\lambda^5}{3125}
\right)  & 
\frac{3\lambda^4}{25000} \;_{2}F_1\left( \begin{array}{c} \frac{6}{5} \; \frac{8}{5} \\ \frac{9}{5} \end{array} ; \frac{-\lambda^5}{3125}
\right)  \\
-{2\lambda  } \;_{2}F_1\left( \begin{array}{c} \frac{2}{5} \; \frac{4}{5} \\ \frac{6}{5} \end{array} ; \frac{-\lambda^5}{3125}
\right)  & 
 \;_{2}F_1\left( \begin{array}{c} \frac{2}{5} \; \frac{4}{5} \\ \frac{1}{5} \end{array} ; \frac{-\lambda^5}{3125}
\right)  
\\ \end{array}\right). \]

\end{Exa}

\begin{Exa}
The final example is the family $X^4+Y^4+Z^4+W^4+\lambda XYZW$. This is a family of $K3$-surfaces. This family is also studied in \cite[pp.  73--77]{DwoPad}.

Considered over a number field, every smooth member of this family has geometric Picard number 19 or 20. This implies that when we consider this family over   a finite field, then every smooth member has geometric Picard number at least 20. From the Tate conjecture (which is proven in this case \cite{NygOg} if $p\geq 5$) it follows that every smooth member has Picard number 20 or 22. This implies that at least 19 of the eigenvalues of $\Frob^*_{q,\lambda}$ on $H^3(U_\lambda)$  are of the form $q\zeta$, with $\zeta$ a root of unity. We will indicate how one can obtain this result from the methods described in this section.

First we calculate the operator $A(\lambda)$. We obtain that
\[ A(\lambda)\omega_{(1,2,2,3)} =  \;_1F_0\left( \begin{matrix} \frac{1}{2} \\ - \end{matrix} ; \frac{\lambda^4}{256} \right) \omega_{(1,2,2,3)}.\]
The operator $A(\lambda)$ leaves the space spanned by $\omega_{(1,1,3,3)}$ and $\omega_{(3,3,1,1)}$ invariant. Its action is as follows:
\[ \left(\begin{array}{cc} 
 \;_{2}F_1\left( \begin{array}{c} \frac{1}{4} \; \frac{3}{4} \\ \frac{1}{2} \end{array} ; \frac{\lambda^4}{256}
\right)  & 
\frac{\lambda^2}{32} \;_{2}F_1\left( \begin{array}{c} \frac{3}{4} \; \frac{5}{4} \\ \frac{3}{2} \end{array} ; \frac{\lambda^4}{256}
\right)  \\
\frac{\lambda^2}{32} \;_{2}F_1\left( \begin{array}{c} \frac{3}{4} \; \frac{5}{4} \\ \frac{3}{2} \end{array} ; \frac{\lambda^4}{256}
\right)  & 
 \;_{2}F_1\left( \begin{array}{c} \frac{1}{4} \; \frac{3}{4} \\ \frac{1}{2} \end{array} ; \frac{\lambda^4}{192}
\right)  
\\ \end{array}\right). \]

One easily computes that
\[   \;_{2}F_1\left( \begin{array}{c} \frac{1}{4} \; \frac{3}{4} \\ \frac{1}{2} \end{array} ; \frac{\lambda^4}{256}
\right)  \pm
\frac{\lambda^2}{32} \;_{2}F_1\left( \begin{array}{c} \frac{3}{4} \; \frac{5}{4} \\ \frac{3}{2} \end{array} ; \frac{\lambda^4}{256}
\right)= \;_1F_0\left( \begin{matrix} \frac{1}{2} \\ - \end{matrix} ; \frac{\pm \lambda^2}{16} \right)\]
hence
\[  A(\lambda)\omega_{(1,1,3,3)}\pm \omega_{(3,3,1,1)} =  \;_1F_0\left( \begin{matrix} \frac{1}{2} \\ - \end{matrix} ; \frac{\pm \lambda^2}{16} \right) \omega_{(1,1,3,3)}\pm \omega_{(3,3,1,1)}.\]
As explained in the previous example, this implies that if $q\equiv 1 \bmod 4$ then $\Frob^*_{\lambda,q}$ restricted to the subspace generated by the  $\omega_{(1,2,2,3)}$, $\omega_{(1,1,3,3)}$ and all the coordinate permutations of these forms, is a (quartic) twist of $\Frob^*_{0,q}$. Using Jacobi sums one can show that the $\Frob^*_{0,q}$ restricted to this subspace has only eigenvalues of the form $q\zeta$, with $\zeta$ a root of unity. This yields 18 eigenvalues of $\Frob_{\lambda,q}$ of this form. Since the number of eigenvalues of $\Frob_{\lambda,q}^*$ that are not of this form  is even, and the complementary subspace has dimension 3, there is a 19-th eigenvalue of the form $q\zeta$.

The final subspace under consideration is $\spa \{\omega_{(1,1,1,1)}, \omega_{(2,2,2,2)}, \omega_{(3,3,3,3)}\}$. We obtain the following matrix with respect to this basis:
\[ \left(\begin{array}{ccc} 
 \;_{3}F_2\left( \begin{array}{c} \frac{1}{4} \; \frac{1}{4} \; \frac{1}{4}\\ \frac{1}{2} \;\; \frac{3}{4}\end{array} ; \frac{\lambda^4}{256}
\right)  & 
\frac{-\lambda^3}{1536} 
 \;_{3}F_2\left( \begin{array}{c} \frac{5}{4} \; \frac{5}{4} \; \frac{5}{4}\\ \frac{3}{2} \;\; \frac{7}{4}\end{array} ; \frac{\lambda^4}{256}\right) & 
 
\frac{\lambda^2}{1024} 
 \;_{3}F_2\left( \begin{array}{c} \frac{5}{4} \; \frac{5}{4} \; \frac{5}{4}\\ \frac{3}{4} \;\; \frac{3}{2}\end{array} ; \frac{\lambda^4}{256}\right) 
 \\
 -\lambda \;_{3}F_2\left( \begin{array}{c} \frac{1}{2} \; \frac{1}{2} \; \frac{1}{2}\\ \frac{3}{4} \;\; \frac{5}{4}\end{array} ; \frac{\lambda^4}{256}
\right)  & 

 \;_{3}F_2\left( \begin{array}{c} \frac{1}{2} \; \frac{1}{2} \; \frac{1}{2}\\ \frac{1}{4} \;\; \frac{3}{4}\end{array} ; \frac{\lambda^4}{256}\right) & 
 
\frac{-\lambda^3}{192} 
 \;_{3}F_2\left( \begin{array}{c} \frac{3}{2} \; \frac{3}{2} \; \frac{3}{2}\\ \frac{5}{4} \;\; \frac{7}{4}\end{array} ; \frac{\lambda^4}{256}\right) 
 \\
\lambda^2 
 \;_{3}F_2\left( \begin{array}{c} \frac{3}{4} \; \frac{3}{4} \; \frac{3}{4}\\ \frac{5}{4} \;\; \frac{3}{2}\end{array} ; \frac{\lambda^4}{256}
\right)  & 
-2\lambda 
 \;_{3}F_2\left( \begin{array}{c} \frac{3}{4} \; \frac{3}{4} \; \frac{3}{4}\\ \frac{1}{2} \;\; \frac{5}{4}\end{array} ; \frac{\lambda^4}{256}\right) & 
 
 \;_{3}F_2\left( \begin{array}{c} \frac{3}{4} \; \frac{3}{4} \; \frac{3}{4}\\ \frac{1}{4} \;\; \frac{1}{2}\end{array} ; \frac{\lambda^4}{256}\right) 
 \end{array}\right) \]
\end{Exa}
\section{Fermat hypersurfaces \& Equivalence relations}\label{SecDia}
In the previous sections it is shown how to calculate the deformation matrix $A(\lambda)$. In this section we discuss the Frobenius action on the  central fiber.

\begin{Lem} \label{LemDia}Let $\mathbf k$ be an admissible monomial type. Let $\mathbf m = \overline{q}\mathbf k$. We have $\Frob_{0,q} \omega_{\mathbf k} = c_{\mathbf k,q} \omega_{\mathbf m}$ for some $c_{\mathbf k,q}$.\end{Lem}
\begin{proof}
Take as a lift of Frobenius the morphism $x_i\mapsto x_i^q$. Then
\[ \Frob^*_{q,0}(\omega)=\frac{x_i^{qk_i+q-1}}{F(x_i^q)^t}\Omega=\sum_{j=0}  c_{t,j} \frac{x_i^{q k_i}  (F^q-F(x_i^q))^j}{F^{qj+t} }.\]
One can easily show that any exponent of $x_i$ in this sum is congruent to $qk_i+q-1 \bmod d_i$. Hence there is only one monomial type $\mathbf{m}$ occurring in the reduction, namely $\overline{q}\mathbf{k}$.
\end{proof}
 \begin{Rem}Suppose $q\equiv 1 \bmod d$. It is well-known that the eigenvalues of Frobenius on $H^n(U)$ are of the form $q^{n-1}/J_{\mathbf k,q}$, where $J_{\mathbf k,q}$ is a  so-called Jacobi sum. Note that the assumption on $q$ implies that $\overline{q}\mathbf k=\mathbf k$. So the set of Jacobi sums coincides with the set of $c_{\mathbf k,q}$ (cf. the Introduction). A stronger result will be proved in the sequel.
\end{Rem}

\begin{Def}Two monomial types are called \emph{strongly equivalent} if and only if their difference is a multiple of the deformation vector. Two monomial types are called \emph{weakly equivalent} if and only if  there exists non-zero multiples of both monomial types that differ by the deformation vector.  
\end{Def}
 
 The characteristic polynomial of Frobenius on the cohomology can be factorized in factors corresponding to the weak-equivalence classes of monomial types:
 \begin{Thm}\label{PrpFactor}
 Let $\mathbf k$ be an admissible monomial type. Let $S$ be the set of monomial types that are weakly equivalent to $\mathbf k$ and let $S'$ be the set of monomial types that are strongly equivalent to $\mathbf k$. Then 
 \[\Frob_{\lambda,q}\omega_{\mathbf k} = \sum_{\mathbf m \in S} c'_{\mathbf m, q} \omega_{\mathbf m}\]
  for some $c'_{\mathbf m, q}\in \Q_q$. In particular, the characteristic polynomial $P(T)$ of Frobenius  on $H^n(U)$ can be factored  as $P(T)=\prod_{[\kv]} P_{[\kv]}(T)$, where the product is taken over all weak-equivalence classes, and $P_{[\kv]}(T)$ is an element of $\Q_q[T]$ of degree equal to the number of distinct admissible monomial types in the weak-equivalence class $[\kv]$.
  
If, moreover, $q\equiv 1 \bmod d$ then  
  \[\Frob_{\lambda,q}\omega_{\mathbf k} = \sum_{\mathbf m \in S'} c'_{\mathbf m, q} \omega_{\mathbf m}\] 
  for some $c'_{\mathbf m, q}\in \Q_q$. In particular, the characteristic polynomial $P(T)$ of Frobenius  on $H^n(U)$ can be factored  as $P(T)=\prod_{[\kv]} P_{[\kv]}(T)$, where the product is taken over all strong-equivalence classes, and $P_{[\kv]}(T)$ is an element of $\Q_q[T]$ of degree equal to the number of distinct admissible monomial types in the strong-equivalence class $[\kv]$.
   \end{Thm}
 \begin{proof}
 Since $\Frob_{\lambda,q}=A(\lambda)^{-1}\Frob_{\lambda,0}A(\lambda^q)$, it suffices to prove that all these three operators leave the subspace $\spa_{m \in S}(\omega_{\mathbf m})$ (if $q \not\equiv 1 \bmod d$) or the subspace $\spa_{m \in S'}(\omega_{\mathbf m})$ (if $q\equiv 1 \bmod d$) invariant.  For $A(\lambda)^{-1}$ and $A(\lambda^q)$ this follows from Proposition~\ref{PrpCoeff}. For $\Frob_{\lambda,0}$ this follows from Lemma~\ref{LemDia}.
 \end{proof}
 \begin{Rem}
 In Corollary~\ref{CorFac} we show that the factorization mentioned above gives factors which are polynomials with $\Q$-coefficients rather then with $\Q_q$-co\-ef\-fi\-ci\-ents. 
 \end{Rem}
 
 It remains to show that `weak equivalence' is the same relation as `nondistinguishable by automorphisms'. 

 \begin{Def}
We call $\mathbf{b}\in (\Z/d\Z)^{n+1}$ an {\em admissible automorphism type} if  $\bv$ can be written as $(w_0b_0,w_1b_1,\dots,w_nb_n) \in(\Z/d\Z)^{n+1}$, such that $\sum w_ib_ia_i \equiv 0 \bmod d$.  Define $\sigma_{\mathbf{b}}$ to be the automorphism
\[ [x_0:x_1:\dots:x_n]\mapsto [\zeta_{d}^{w_0b_0} x_0:\zeta_{d}^{w_1b_1}x_1: \dots:\zeta_d^{w_nb_n}x_n].\]
We call two monomial types $\mathbf{k}$ and $\mathbf{m}$ {\em distinguishable by automorphisms} if there exists an admissible automorphism type $\mathbf{b}\in (\Z/d\Z)^{n+1}$ such that
\[  \sigma_{\mathbf{b}}\left(\prod x_i^{k_i}\right) =\prod x_i^{k_i} \mbox{ and } \sigma_{\mathbf{b}}\left(\prod x_i^{m_i} \right) \neq \prod x_i^{m_i}. \]
\end{Def}

\begin{Thm} \label{ThmWeak}Two monomial types $\mathbf{k}$ and $\mathbf{m}$ are weakly equivalent if and only if $\mathbf{k}$ and $\mathbf{m}$ are not distinguishable by automorphisms.\end{Thm}

 \begin{proof} One easily sees that $\sigma_{(b _i)}$ fixes $\omega_\kv$ if and only if $\prod x_i^{k_i+1}$ is fixed by
$\sigma_{({b_i})}$. This in turn is equivalent with 
 \[ \sum b_i (k_i+1) w_i \equiv 0 \bmod d,\]
  and similarly for $\mathbf{m}$.

 (`$\Rightarrow$'). Suppose $\kv$ and $\mv$ are weakly equivalent. Then we have a relation
 \[ s \mathbf{k} + t \mathbf{m} =r \mathbf{a},\] with $s,t\in (\Z/d\Z)^*$.
 It suffices to show that if $\bv$ is an admissible automorphism type then 
 \[ \sum b_i (k_i+1) w_i \equiv 0 \bmod d \Leftrightarrow \sum b_i (m_i+1) w_i \equiv 0 \bmod d.\]
  Since $\kv$ and $\mv$ are weakly equivalent we have
\[ s \sum  b_iw_i (k_i+1) + t \sum  b_i w_i (m_i+1) \equiv r \sum b_i a_i w_i \equiv 0 \bmod d.\]
Hence
 \[ s \sum b_i w_i (m_i+1)  \equiv  -t \sum  b_iw_i (k_i+1) \bmod d.\] Since $s$ and $t$ are invertible, the above claim follows.
 
(`$\Leftarrow$'). Suppose $\kv$ and $\mv$ are not distinguishable by automorphisms. Take $\bv$ such that $\sigma_\bv(\omega_\kv)=\omega_\kv$ and $\sigma_\bv(\omega_\mv)\neq \omega_\mv$
 Hence 
 \[ \sum b_i (k_i+1) w_i \equiv 0 \bmod d,\]
 and
 \[ \sum b_i (m_i+1)w_i \not \equiv 0 \bmod d.\]
Suppose $\kv$ and $\mv$ are weakly equivalent, i.e., we have a relation 
\[ s \mathbf{k} + t \mathbf{m} =r \mathbf{a}\]
where $s$ and $t$ are invertible in $\Z/d\Z$.
Then 
\[ s \sum  b_iw_i (k_i+1) + t \sum  b_i q_i (m_i+1) - r \sum b_i a_i q_i \equiv 0 \bmod d.\]
Since the first and third summand are zero, the same holds for the second summand. Contradicting that it should be non-zero. So we cannot have a relation
\[ s \kv + t \mv =r \av.\]
Hence $\mathbf{k}$ and $\mathbf{m}$ are not weakly equivalent.
 \end{proof}

\begin{Def} Assume that $q\equiv 1 \bmod d$ (i.e., $\F_q\supset \F_p(\zeta_d)$). Let $\chi$ be the $d$-th power residue symbol. Let $\kv$ be an admissible monomial type. Let $\kv_i$ be the $i$-th entry of $\kv$, i.e., $w_i(k_i+1)$. Then the Jacobi-sum associated with $\kv$ is defined as
\[ J_{\kv,q}:=(-1)^{n+1} \sum_{(v_1,\dots,v_n)\in \F_q^n \colon \sum_i v_i=-1} \chi(v_1)^{\kv_1} \chi(v_2)^{\kv_2}\dots \chi(v_n)^{\kv_n}.\]
\end{Def}
 
\begin{Cor}\label{CorJac} Assume $q\equiv 1 \bmod d$. Let $\mathbf k$ be an admissible monomial type. Let $S$ be the set of monomial types that cannot be distinguished by automorphisms from $\mathbf k$. Then the sets $S_1:=\{ q^n/c_{\mathbf m,q} \colon m \in S\}$ and $S_2:=\{J_{\mathbf m, q}  \colon m \in S\}$ coincide.

\end{Cor} 
\begin{proof}
Let $G\subset \prod \Z/d_i\Z$ be the group of automorphisms that fixes $\omega_{\mathbf k}$. Then $X_0/G$ is  a Fermat variety in a different weighted projective space $\Ps'$. It is well-known that the eigenvalues of Frobenius on the primitive part of $H^{n-1}_{\rig,c}(X_0/G)$ are Jacobi-sums appearing in $S_2$ (e.g. see \cite{GY}).

The group $H^n(U/G)$ is canonically isomorphic with the subspace of $H^n(U)$ generated by the forms $\omega_{\mathbf m}$, where $\mathbf m \in S$ (this follows from \cite[Lemma 2.2.2]{Dol}). This implies that all the $q^n/c_{\mathbf m,q}$  with $\mathbf m \in S$ are  eigenvalues of $\Frob$ on $H^{n-1}_{\rig,c}(X_0/G)$. Hence
$S_1=S_2$.
\end{proof}
 
\begin{Cor}
\label{CorFac} Let $\overline{\lambda} \in \F_q$.  Let $P(t)$ be the characteristic polynomial of $\Frob_{\overline{\lambda}}$ on $H^n(U_\lambda,\Q_q)$. Then 
\[ P(t)=\prod_{[\kv]}P_{[\kv]}(t), \]
where the product is taken over all weak-equivalence classes of admissible monomial types. Let $\kv$ be an admissible polynomial type. Then $P_{[\kv]}(t)$ is an element of $\Q[t]$ and its degree equals the number of admissible monomial types that are weakly equivalent with $\kv$.  
\end{Cor}

\begin{proof}
Fix for the moment a monomial type $\kv$. Let $G_\kv \subset \prod \Z/d_i\Z$ be the group of automorphisms that fixes $\omega_{\mathbf k}$. Then $X_0/G_\kv$ is  a Fermat variety in a different weighted projective space $\Ps'$ and $H^n(U_0/G_\kv,\Q_q)$ is canonically isomorphic with the subspace of $H^n(U_0,\Q_q)$ generated by the form $\omega_{\mathbf m}$, where $\mv$ is weakly equivalent with $\kv$.
This enables us to write 
\[ H^n(U) = \oplus_{[\kv]} H^n(U/G_{[\kv]}). \]
For every weak-equivalence class of monomial types, set  $P_{[\kv]}(t)\in \Q[t]$ to be the characteristic polynomial of Frobenius acting on $H^n(U/G_{[\kv]})$. Then $P(t)=\prod P_{[\kv]}(t)$, we have
$P_{[\kv]}(t)\in \Q[t]$ and
 \[\deg(P_{[\kv]}(t))=\dim H^n(U/G_{[\kv]}) = \# \{ \mv \colon \kv \mbox{ and } \mv \mbox{ are weakly equivalent} \}, \]
 which finishes the proof.
\end{proof}

 \begin{Exa}\label{exaQuiMon} Consider the case of the quintic threefold in $\Ps^4$, with deformation vector $a=(1,1,1,1,1)$.
 Up to interchanging coordinates  we have the following five strong equivalence classes:
 \begin{enumerate}
 \item $[0,0,0,0,0], [1,1,1,1,1], [2,2,2,2,2], [3,3,3,3,3]$.
 \item $[0,0,0,2,3], [3,3,3,0,1]$. (20 Permutations possible)
 \item $[1,1,1,0,2], [2,2,2,1,3]$. (20 Permutations possible)
 \item $[0,0,1,1,3], [2,2,3,3,0]$. (30 Permutations possible)
 \item $[2,2,0,0,1], [3,3,1,1,2]$. (30 Permutations possible)
 \end{enumerate}
 The classes (2) and (3) form one weak-equivalence class, the same holds for (4) and (5). 
 Over an arbitrary finite field we obtain three distinguishable factors of the zeta-function, all three of degree 4. One factor is occurring with multiplicity 30, one factor is occurring  with multiplicity 20, and one factor is occurring with multiplicity one. This is in agreement with \cite{Cand2}. 
 \end{Exa}

\end{document}